\documentclass[12pt]{amsbook}
\usepackage{amsmath}

\theoremstyle{plain}
\newtheorem{thm}{Theorem}[chapter]
\newtheorem{lemma}[thm]{Lemma}
\newtheorem{cor}[thm]{Corollary}
\newtheorem{prop}[thm]{Proposition}

\theoremstyle{remark}

\def \TH {\Theta}
\def \tcG {\tilde{\mathcal{G}}}
\def \tcB {\tilde{\mathcal{B}}}
\def \tg {\tilde{g}}
\def \sh {\sharp}
\def \Dad {D_{A^\a}^\de}
\def \dap {d_{A^\a}^{+}}
\def \datp {d_{A^\a}^{\tilde{+}}}
\def \tdad {d_{A^\a}^{\tst_\de}}

\def \tmu {\tilde{\mu}}
\def \tm {\tilde{-}}
\def \tp {\tilde{+}}
\def \tM {\tilde{M}}
\def \NK {N \! \setminus \! K}
\def \ad {\text{ad}}
\def \cH {\mathcal{H}}
\def \tcH {\tilde{\cH}}
\def \su  {\mathfrak{su}_2}
\def \lra {\longrightarrow}
\def \cpi {\frac{1}{8\pi^2}}
\def \tb {\tilde{b}}
\def \wCS {\widetilde{CS}^\a}
\def \deg {\text{deg}}
\def \t {\tau}
\def \ker {\text{ker}}
\def \coker {\text{coker}}
\def \im {\text{im}}

\def \w {\wedge}
\def \aa{\accent'27a}
\def \Aa {\mathcal{A}^\alpha}
\def \tX {\tilde{X}}
\def \e {\epsilon}

\def \Om {\Omega}
\def \L {\Lambda}

\def \k {\kappa}
\def \ra {\rightarrow}
\def \r {\rho}
\def \de {\delta}
\def \sb {spec_b}
\def \rsb {\Re \sb}

\def \z {\zeta}
\def \R {\mathbb{R}}
\def \Z {\mathbb{Z}}
\def \C {\mathbb{C}}
\def \Sig{\Sigma}
\def \S {\Sigma}
\def \d {\partial}
\def \a {\alpha}
\def \LL {L^{\otimes 2}}
\def \D {\Delta}
\def \b {\beta}
\def \s {\sigma}
\def \x {\xi}
\def \g {\gamma}
\def \h {\eta}
\def \L {\Lambda}
\def \Vedge {\mathcal{V}_e}
\def \n {\nu}
\def \th {\theta}

\def \cg {\check{g}}
\def \fg {\mathfrak{g}}
\def \XS {X \!  \setminus \! \Sig}
\def \Lpka {L^p_{k, A^\a}}
\def \Lpca {L^p_{1, A^\a}}
\def \Lpga {L^p_{2, A^\a}}

\def \cDa {\check{\D}_{\a,0}}
\def \tHgn {\tilde{H}_{g,n}}
\def \tth {\tilde{\th}}
\def \cst {\check{*}}

\def \ak {\tfrac{\a}{\k}}
\def \ka {\tfrac{\k}{\a}}
\def \rr {\frac{1}{r}}
\def \da { ( \d_\th + 2i \k ) }
\def \rcDa {r^2 \cDa}
\def \ma {\tfrac{m}{\a}}
\def \l {\lambda}
\def \half {\tfrac{1}{2}}
\def \tDa {\tilde{\D}_\a}
\def \tst {\tilde{*}}
\def \tda {d_{A^\a}^{\tst}}
\def \G {\Gamma}
\def \aD {a_{\! D}}
\def \tD {\tilde{\D}}
\def \td {d^{\tst}}
\def \cDaf {\check{\D}_{\a,1}}
\def \cDak {\check{\D}_{\a,k}}

\def \tDak {\tilde{\D}_{\a,k}}
\def \rcDaf {r^2 \cDaf}
\def \pDa {\D'_{\a,1}}
\def \rpDa {r^2 \pDa}
\def \Bap {\mathcal{B}^{\a,p}}
\def \tA {\tilde{A}}
\def \tE {\tilde{E}}
\def \tL {\tilde{L}}
\def \tr {\text{tr}}
\def \MK {M \! \setminus \! K}
\def \TK {T \! \setminus \! K}
\def \om {\omega}
\def \cFa {\mathcal{F}^\a}
\def \WS {W \! \setminus \! \S}
\def \MKj {M_j \! \setminus \! K_j}
\def \MjR {M_j \times \R^+}
\def \tDM {\tD^{M}}
\def \tDW {\tD^{W}}
\def \tL {\tilde{L}}
\def \hD {\hat{\D}}

\begin{document}
\frontmatter
\title{Singular Connections on Three-Manifolds and Manifolds with
Cylindrical Ends.}
\author{Benoit G\'{e}rard\footnote{Email address: Benoit\_Gerard@MCKINSEY.COM}}
\maketitle

\tableofcontents

\mainmatter

\chapter*{Introduction.}

This thesis can be seen as a first step towards defining a Floer homology
for the singular connections studied by Kronheimer and Mrowka in their
work on the structure of Donaldson's polynomial invariants \cite{KM}.  By
investigating certain spaces of connections on a homology three-sphere which
are singular along a knot, such a Floer homology would provide possibly a new
invariant for the knot.  We will first briefly recall the key ideas of the
classical instanton homology defined by Floer \cite{Floer} for a homology 
three--sphere and describe the singular connections considered by Kronheimer
and Mrowka.  We will then present the content of the thesis.

Floer homology is obtained by applying an infinite--dimensional version of
Morse theory to the Chern-Simons functional defined over the space of
connections in a principal bundle (with e.g. structure group $SU(2)$) over
a closed, oriented, three--manifold $M$.  The set $\mathcal{R}$ of critical
points of this functional are the flat connections and the gradient-flow
lines correspond to anti--self--dual connections over the cylinder $M \times \R$.
After generic perturbation, the critical points form a finite set of isolated
connections.  The index of the critical points is always infinite, but it is
nevertheless possible to define a relative index between two critical
points by the spectral flow \cite{APS} of the Hessian of the Chern-Simons
functional.  After dividing by the action of the gauge group, this spectral flow
is well-defined modulo $8$.  It follows that the free abelian  group $R$
generated by the set of critical points is $\Z / 8$--graded.  If we consider
two critical points $\a$ and $\b$ with consecutive gradings, the space of flow 
lines starting from $\a$ and arriving at $\b$ is generically a zero--dimensional,
oriented, compact manifold, i.e. a finite set of signed points.  Counting those
flow lines as in classical finite--dimensional Morse theory yields a boundary
map $\de : R \ra R$, and Floer proved that $\de^2 = 0$ by analyzing the ends of the
spaces of flow lines.  Thus, one gets an invariant for the manifold $M$ in the form
of $8$ homology groups $HF_*(M) = \ker \, \de / \im \, \de$.  Floer homology
is a basic tool in computations of Donaldson's polynomial invariants \cite{DK}
for four--manifolds , see e.g. \cite{MMR}, \cite{TaubesL2}.

On the other hand, Kronheimer and Mrowka \cite{KM} developed a gauge theory
for a surface $\S$ embedded in a four--manifold $X$.  They study spaces of
ASD connections on $X$ singular along $\S$, but whose holonomy around small
loops linking $\S$ is asymptotically fixed by a parameter $\a \in (0, \half)$.
They then define a numerical invariant similar to Donaldson's polynomials.
They use that invariant to establish certain constraints on the genus of $\S$, which
are reminiscent of the adjunction formula for curves in complex surfaces.

The first chaper of this thesis starts by recalling the set-up introduced in
\cite{KM} and some results obtained there regarding the moduli spaces of 
singular ASD connections.  The analysis in \cite{KM} is based on orbifold
techniques.  They put an orbifold metric with a sharp cone--angle 
$\frac{2 \pi}{\nu}, \nu \in \Z^+$, around $\S$, and start from a
rational holonomy parameter $\a = \frac{p}{\nu}, p \in \Z^+$.  They then
use a trick from \cite{DS} to change the metric back to the round
metric.  An estimate based on a Weitzenbock formula is used to change the
holonomy parameter to any real value of $\a$ within a certain compact
interval of $(0, \half)$ determined by $\nu$.  Instead of orbifold techniques,
we use the theory of elliptic edge operators developed in \cite{Mazzeo} and
summarized in the appendix  to this thesis.  Edge theory has the advantage
that one can work directly with any sharp enough cone-angle and any real value
of the holonomy parameter $\a$.  This in turn allows one to write down an
explicit global Coulomb gauge--fixing condition, something which doesn't appear
in \cite{KM}.  Another problem with the orbifold techniques is that the estimate
which allows one to go from a rational value to any real value of $\a$ is  based
on the fact that the four--sphere, used as a local model, has positive Ricci 
curvature.  That positive Ricci curvature assumption cannot be made when one
works with manifolds with cylindrical ends.  Hence, the perturbation
argument for $\a$ does not carry over and this is why one really needs another
analytical theory.  The last section of chapter $1$ carries out some more
edge analysis which will be used in the subsequent chapters.

In the second chapter, we first show how to adapt the setting of \cite{KM}
to study singular three--dimensional connections.  We then show how to define
a Chern--Simons functional $CS^\a$ for singular connections.  As usual, the
critical points of the Chern--Simons functional are flat connections and
the gradient--flow lines are singular ASD connections.  The last section
of chapter $2$ computes the kernel of the Hessian of $CS^\a$ at a critical
point, which yields a topological condition for determining when a
singular flat connection is isolated.  This computation is in
fact the Hodge--de Rham theory of singular flat connections.  The chapter ends
with a word about perturbations of the Chern--Simons functional.

Chapter $3$ treats the case of singular connections on manifolds with cylindrical
ends.  After establishing the set--up, we study the local models on the ends.
The analysis is reminiscent of \cite{LM} and relies on the results 
obtained at the end of chapter $1$.  It is also here that the trick of 
changing the conformal strucure from \cite{DS} is used.  The global theory
of moduli spaces of ASD singular connections is then presented.  The chapter 
ends with a few index formulae for gluing and changing gauge, which
show that an analogue of Floer's grading function can be defined
for singular flat connections.  However, in this case, the grading function
is $4$--periodic, not $8$--periodic.

The thesis ends with a few comments regarding the work that remains to be done
to define a Floer homology and a couple of conjectures regarding possible 
relations between that putative Floer homology and other gauge theories for knots.

There is an appendix which contains a summary of the theory of elliptic edge
operators.  Moreover, the theory is there extended to $L^p$--spaces and a number
of results regarding elliptic estimates for the operators arising in the thesis
are presented.

\vfil\eject
\chapter{Four-dimensional Theory.}

This chapter deals with the theory of anti-self-dual connections on
a $4$--manifold which are singular along an embedded surface.  The first section
recalls the set-up developed by Kronheimer and Mrowka \cite{KM} to study the
moduli spaces of such connections and states some of the results they
obtained (dimension of the moduli space, Chern--Weil formula, etc.).
The second section develops a slice theorem for the action of the gauge
group on the space of singular connections. The analytical techniques
used in that section involve the theory of elliptic edge operators 
developed in \cite{Mazzeo}.  The reader who is not knowledgeable about
edge theory will find a summary in the appendix to this thesis.
The last section is an analysis of certain twisted Laplacians
corresponding to singular connections on differential forms.  The results
obtained there will be useful in chapter $2$ to study the Hodge--de Rham
theory of flat singular connections on $3$--manifolds and in chapter $3$
to develop the moduli theory of singular connections on manifolds with
cylindrical ends.

\section{Set-up.} \label{S:setup}

This section recalls a number of facts from
Kronheimer and Mrowka's article \cite{KM} for $SU(2)$ connections on a 4-manifold
$X$, singular along a surface $\Sig$.  It recalls the set-up they developed and some of 
the results they obtained regarding the moduli spaces of singular anti-self-dual connections.

\subsection{Basic definitions.} \label{SS:fourdimsetup}

 Let $X$ be a smooth, closed, oriented 4-manifold, endowed with a Riemannian
metric.  Let $\Sig$ be a closed, oriented embedded surface of genus $g$ 
and $\n$ the 
normal bundle of $\Sig$ in X.  Using the exponential map given by the metric
on $X$, we identify a small disk bundle of $\n$ with a closed tubular neighborhood
$N$ of $\Sig$.  Let $r$ be the variable defined by the distance from the 
zero section in $\n$.

 Let $E$ be a $SU(2)$-bundle over $X$; we take an abelian
reduction of E over $N$
\begin{equation}  \label{E:splitting}
       E|_N = L \oplus L^{-1}, 
\end{equation}
where $L$ is a $U(1)$-bundle on $N$.

Define the two topological numbers
$$
 k \: = \: c_2(E)[X],
$$
$$
  l \: = \: - c_1(L)[\Sig].
$$

  Let $A^0$ be a connection on E, which
respects the splitting \eqref{E:splitting}, i.e. which is abelian on $N$. 
Now, $\d N$ is a circle bundle over $\Sig$, for which we choose
a connection 1-form $i\h$.  Let $\pi : N \setminus \Sig \, \ra \, \d N$ be the radial
projection.  Then $\pi^{*} \h$, which, by abuse of notation, we will also denote
$\h$, is an angular 1-form on $N$.  As in \cite{KM},
 define a ``background'' connection on the 
restriction of $E$ to $\XS$ by
\begin{equation}   \label{E:backconn}
  A^{\a} \: = \: A^0 + i \b (r)
\begin{pmatrix} 
 \a & 0 \\ 0 & -\a   
\end{pmatrix}
\h,
\end{equation}
where $\a \in (0, \frac{1}{2})$ and 
$\b$ is a cut-off function equal to $1$ near $0$.  By construction,
$A^\a$ also respects the splitting \eqref{E:splitting} over $N$. The holonomy of $A^\a$
around a small loop linking $\Sig$ tends to
\begin{equation}  \label{E:holonomy}
\begin{pmatrix}
  e^{- \! 2i\pi \a}  & 0  \\
  0 & e^{2i\pi \a}
\end{pmatrix}
\end{equation}
as the loop shrinks towards $\Sig$.  

Let $\fg_E$ denote the adjoint bundle of $E$.  We denote by $\Lpka$ the Sobolev spaces
defined with
the adjoint connection of $A^\a$ on $\XS$ on sections
 or forms with values in $\text{End}\, E$ or $\fg_E$.

 For $p > 2$, define the space of singular connections
\begin{equation}  \label{E:spaceconn}
 \mathcal{A}^{\a,p} (X, \Sig) \: = \: A^\a + \Lpca (\Om^1( \XS \, ; \fg_E))
\end{equation}
and the corresponding space of gauge transformations
\begin{equation}  \label{E:spacegauge}
  \mathcal{G}^p (X, \Sig) \: = \: \{ g \in \text{Aut}(E) \, \mid \, g \in \Lpga( \XS \, ; \text{End} \, E) \}.
\end{equation}
For reasons that will be recalled later, $\mathcal{G}^p$ is independent of $\a$.

The space $\mathcal{G}^p$ is a Banach Lie group,
 acting smoothly on $\mathcal{A}^{\a,p}$.
 The Lie algebra of $\mathcal{G}^p (X, \S)$ is $L^p_{2,A^\a}( \XS, \fg_E )$.
The resulting quotient is denoted
$\mathcal{B}^{\a,p}(X, \Sig) = \mathcal{A}^{\a,p} / \mathcal{G}^p$. The subspace $( \mathcal{B}^{\a,p} )^*$
of $\mathcal{B}^{\a,p}$ consisting of irreducible connections  is a smooth Banach manifold.
We will give a new proof of this fact, by
constructing an explicit slice for the action of the gauge group, see section 
\ref{S:slice}.

We define
$$
  M^{\a,p}_{k,l} (X, \Sig) = \{ [A] \in \mathcal{B}^{\a,p}(X, \S) \, \mid \, F^+_A = 0 \} \, ,
$$
to be the  moduli space of singular anti-self-dual connections.
The deformation complex for the anti-self-duality
equation is 
\begin{equation}   \label{E:ASDcomplex}
 0 \ra L^p_{2,A^\a} (\XS ; \fg_E) \overset{d_A}{\lra} L^p_{1,A^\a} ( \Om^1( \XS ; \fg_E))
   \overset{d_A^+}{\lra} L^p(\Om^+ (\XS ; \fg_E)) \ra 0.
\end{equation}
Kronheimer and Mrowka \cite{KM} show that 
the index of this complex, that is the formal dimension of the
moduli space  $M^{\a,p}_{k,l}$, is given by
\begin{equation}   \label{E:dimASD}
   d_{k,l} = 8k + 4l - 3(b^+_2 (X) - b_1(X) + 1) - (2g - 2).
\end{equation} 
They also show that $M^{\a,p}_{k,l}$ 
is independent of $p$ for $p$ close enough to $2$, and that 
for a generic metric on $X$, this moduli space is cut-out transversely 
away from the flat connections, and is smooth away from the reducible
or flat connections.

We will also need the Chern-Weil formula for singular connections in 
$\mathcal{A}^\a$ \cite[proposition 5.7]{KM} :
\begin{equation}  \label{E:CW}
 \frac{1}{8\pi^2} \int_{ X \setminus \Sig } \text{tr} \, (F_A \wedge F_A ) \: = \:
    k + 2\a l - \a^2 \Sig \cdot \Sig .
\end{equation}

\subsection{Cone-like metrics}   \label{SS:cone}
As in \cite{KM}, the analysis for the $L^{\otimes 2}$ component of the adjoint
bundle $\fg_E$ near the surface $\Sig$ will be carried out by putting a 
cone-like metric around $\Sig$. In \cite{KM}, the authors restrict themselves
to cone-angles which are divisors of $2\pi$,
in order to use orbifold analysis techniques.
Instead, we will consider general real cone-angles, and use the theory
of edge differential operators developed in \cite{Mazzeo}, as presented
and slightly extended in appendix \ref{A:edge}.

The way to construct such a cone-like metric is presented in \cite{KM}.  We
are using the notations introduced in \ref{SS:fourdimsetup}.  Let
$g$ be the metric on $X$, and $g_\Sig$ the metric induced  on $\Sig$.
 Let $p: \nu \ra \Sig$ be the natural projection.  Then,
the metric
$$
  \bar{g}_\e \; = \;  (1-\b (r/ \e)) g \, + \,
        \b (r/ \e) ( p^{*}g_\Sig + dr^2 + r^2 \h^2 ) 
$$
can be made as close as we want to the metric $g$ in the $C^0$-norm, over
$\XS$, by letting $\e \ra 0$. A cone-like metric $\cg_\e$ with cone-angle
$\ak 2\pi$ around $\Sig$, is defined by:
$$
  \cg_\e \; = \;  (1-\b (2r/ \e)) \bar{g}_\e \, + \,
      \b (2r/ \e) \biggl( p^{*}g_\Sig + dr^2 + \bigl( \ak \bigr) ^2 r^2 \h^2 \biggr) .
$$
The norms defined on $L^p$-spaces from the different metrics $g$, 
$\bar{g}_\e$ and $\cg_\e$ are of course different, but since these metrics
are comparable, they will be equivalent.

Those different metrics yield different conformal structures on
$\XS$. As explained in detail in \cite{DS}, there is a simple way to
understand how a change in the conformal structure of a four-manifold
affects the space of anti-self-dual
$2$-forms.  Given a background metric on
a four-manifold, with self-dual space $\L^+$ and anti-self-dual space $\L^-$, the 
anti-self-dual space for another metric will be given by the graph of a unique
linear map
$$
   \mu : \L^- \ra \L^+ ,
$$
the operator norm of which is pointwise everywhere less than $1$.  As in \cite{KM}, it
is straightforward to compute that the change in the
conformal structure from $\cg_\e$ to $\bar{g}_\e$ will be given by a map $\mu$ supported
on an $\frac{\e}{2}$-neighborhood of $\Sig$ and such that
\begin{equation}  \label{E:normmu1}
  | \mu | \leq \frac{\k / \a -1}{\k / \a +1} = \frac{\k - \a}{\k +  \a}
\end{equation}
pointwise.  Since the metrics $\bar{g}_\e$ and $g$ can be made as close as we want, changing
the conformal structure from $\cg_\e$ to $g$ entails a map $\mu$ with a bound only slightly
bigger than the one given by
\eqref{E:normmu1}, i.e.
\begin{equation}  \label{E:normmu}
    | \mu | \leq \frac{\k - \a}{\k+ \a} +  \de \, ,
\end{equation}
where $\de \ra 0$ as $\e \ra 0$.  From now on, we'll drop the subscript $\e$.

It is also possible to change the cone-angle from $\tfrac{\a}{\k}$ to $\tfrac{\a}{\k'} ,
\k' \leq \k$.  In that case, the corresponding $\mu$ map satisfies the bound
\begin{equation}   \label{E:intcone}
  | \mu | \leq \frac{\k - \k'}{\k + \k'}.
\end{equation}

Finally, we record the following fact regarding the Hodge star operator $\cst$ in
the cone-angle metric :
\begin{equation}   \label{E:regconeHodge}
  \cst \in  \bigcap_{\ell = 0}^{ \infty } \frac{1}{r^\ell} L^\infty_\ell 
        ( \XS ; \text{End} \, \L^* ) ,
\end{equation}
where the $L^\infty_\ell$ spaces are defined with the Levi-Civita connection
corresponding to the round metric.

\subsection{Function spaces}   \label{SS:fspaces}

We recall here a number of facts regarding the function spaces $L^p_{k,A^\a}$
introduced in \ref{SS:fourdimsetup}.  Let $F$ be any vector
bundle over $X$, and let $\nabla$ be a connection on $F$.  We extend the distance
function $r$  smoothly to a function that takes the value $1$ on $X \! \setminus \! N$.  Then, as in
\cite{KM}, 
we define $W^p_k(X;F)$ as the completion of the space of smooth sections of $F$, compactly
supported on $\XS$, in the norm
\begin{equation}   \label{E:defWpk}
  \| s \|_{W^p_k} \; = \;
      \left\| \frac{1}{r^k} s \right\|_p + \left\| \frac{1}{r^{k-1}} \nabla s \right\|_p
     + \dots +
       \| \nabla^k s \|_p .
\end{equation}
In fact, \cite[lemma 3.2]{KM} shows that
$$
  W^p_k =  r^k L^p \cap r^{k-1} L^p_1  \cap \dots \cap L^p_k
$$
We need to know what those spaces are in terms of the edge Sobolev spaces introduced
in appendix \ref{A:edge}, to which we refer for notation.  Let $\tX \ra X$ be the blow-up of $X$ along
the submanifold $\Sig$ in the category of real, oriented, smooth manifolds.
By this, we mean that we introduce local polar coordinates around
$\Sig$ and declare them to be non-singular.  Actually, $\tX$ is diffeomorphic to
 $X \! \setminus \text{int} \, N$,
but not in any natural way.  The manifold $\tX$ is compact with 
boundary $\d \tX$.  The boundary $\d \tX$ is naturally a circle bundle
over $\Sig$, and this is the edge structure we will consider.  We  pull-back
the bundle $F$ to $\tX$. Then, we have :

\begin{prop} \label{P:eqedge}
 $W^p_k(X;F)$ is naturally isomorphic to $r^{k-1/p}W^{p,k}_e(\tX;F)$.
\end{prop}
\begin{proof}
Since $X$ and $\tX$ are compact, we can work locally in a product neighborhood $U$ of a small
open subset of $\Sig$. 
 Hence, we can  assume that
the bundle $F$ is trivial and that the metric on $X$ is euclidean. Let $(s,t)$ be local
 coordinates on $\Sig \cap U$ extended to $U$, $(x, y)$ local
coordinates on the directions normal to $\Sig$ and $(r, \th )$ the polar 
coordinates corresponding to $(x, y)$. The volume form on $X$ is given by
 $dVol_X = ds \w dt \w dx \w dy \; = \; ds \w  dt \w dr \w r d \th$. Since  we are working
 locally, it is no restriction to assume that the coordinates $(s,t,r, \th)$ are
 orthonormal on
$\tX$ so that its volume form  is given by
$dVol_{\tX} = ds \w dt \w dr \w d \th \; = \; \frac{1}{r} dVol_X$.  This shows already that
$L^p(X) = r^{-1/p}L^p(\tX)$.

Now, still by compactness, it is no restriction to take the connection $\nabla$
to be trivial over $U$. Hence, in polar coordinates, 
it acts as
$$
  \nabla f = \frac{\d f}{\d s} ds +  \frac{\d f}{\d t}dt + \frac{\d f}{\d r} dr
        + \frac{1}{r}\frac{\d f}{\d \th} r 
d \th.
$$
Since $\{ ds, dt, dr, rd\th \} $ form a local orthonormal frame for $\Om^1(X)$
 and $\{ r\d_s, r\d_t, r\d_r, \d_\th \} $ generate
$\Vedge (\tX)$ locally, this shows that $W^p_1(X) = r^{1-1/p}W^{p,1}_e(\tX)$.  The
proposition for $k$ a non-negative integer follows by induction and the Leibniz rule, noting that
$ds, dt \in C^\infty(U)$ and that 
$$
   dr, \, rd\th \quad \in  \quad \bigcap_{\ell=0}^\infty \frac{1}{r^\ell}L^\infty_\ell(U).
$$
For general real values of $k$, use duality and interpolation.
\end{proof}

Recall that, on a tubular neighborhood $N$ of $\Sig$, we have fixed a splitting
$E = L \oplus L^{-1}$.  Therefore, near $\Sig$, the adjoint bundle $\fg_E$ splits
as $i\R \oplus L^{\otimes 2}$.  We refer to $i\R$ as the  \emph{diagonal} component
of $\fg_E$ near $\Sig$ and to $ L^{\otimes 2}$ as the \emph{off-diagonal} component.
Remark that near $\Sig$, the connection $\nabla_{A^\a}$ is
trivial on the diagonal component.

In \cite[proposition 3.7]{KM}, it is shown that $\Lpka$ consists of sections or
forms that are in $L^p_k(X)$, but whose off-diagonal components lie
in $W^p_k(X)$. Therefore,  the space  $\Aa$ consists of connections 
of the form $A^\a + a$, where $a \in L^p_1( \Om^1(X; \fg_E))$ and
 the off-diagonal component of
$a$ lies in $W^p_1( \Om^1(\LL))$ near $\Sig$.  The space of gauge transformations $\mathcal{G}^p$
 consists of
$L^p_2$-sections of $\text{Aut} (E)$ whose off-diagonal component lies in
$W^p_2(\LL)$ near $\Sig$.  Note that this shows why the gauge group is independent
of the holonomy parameter $\a$.  It also follows that the Lie algebra of  $\mathcal{G}^p$ is 
given by the space of sections $u \in L^p_2(X; \fg_E)$ whose off-diagonal component lies in 
$W^p_2( \LL)$ near $\Sig$. 

From now on, we will denote the spaces $\Lpka$ by $\tL^p_k$, emphasizing that they are
independent of $\a$.  The usual Sobolev multiplication and compact embedding theorems
hold for the $\tL^p_k$ spaces, as explained
in \cite[lemma 3.8]{KM}. One just needs to consider things component by component and apply
the usual theorems for unweighted Sobolev spaces.

\section{A slice theorem for singular connections.}  \label{S:slice}

This section develops a Coulomb gauge conditions for singular connections.
This yields a slice theorem for the action of the gauge group, endowing
$\mathcal{B}^{*}$ with the structure of a smooth Banach manifold.  An important point
here is that the slice is given by the kernel of a differential operator,
in fact the adjoint of the covariant derivative with respect to the
$L^2$-inner product given by the round metric on the diagonal component of $\fg_E$ and
an appropriate cone-angle metric on the off-diagonal component. 
In \cite{KM}, the slice wasn't constructed explicitly.  The fact that
the slice can be so given is very useful for gluing theorems, as we will
illustrate in the case of manifolds with cylindrical ends.  Note one serious
difference between the usual
Coulomb condition for non-singular connections $ ( B \in \ker \, d_A^{*} ) $
and the Coulomb condition that we present: the latter is not gauge invariant,
 due to the use of different metrics for the on- and off-diagonal components.

The slice will, as usual, be obtained from analyzing the twisted Laplacian
on sections of $\fg_E$. We will investigate that Laplacian for the
background connection $A^\a$, but first we need to develop a model
problem for the off-diagonal component of $A^\a$, using the theory
of edge differential operators.  The analysis will
be developed in the 4-dimensional case, but it is entirely similar
in the 3-dimensional case.

Note that R\aa de \cite{R1} obtained another gauge fixing condition using
a metric asymptotically isometric to $\mathbb{H}^3 \times S^1$ to analyze the 
off--diagonal component of the connection.  However, his technique fails to
provide a well--defined Coulomb condition when the connection is not abelian
near the embedded surface.

\subsection{The model problem for the off-diagonal component}  \label{SS:modeloff}
 
Before describing the model for the off-diagonal component near
$\Sig$, recall that we have a surface $\Sig$ of genus $g$, embedded in a
four-manifold $X$ with self-intersection $n=\Sig^2$.  Take as a local
model the ruled surface $H_{g,n}$ , obtained by
compactifying the complex line bundle of degree $-n$ over the surface  $\Sig$
with a copy of $\S$ at infinity.
Fix a diffeomorphism between the tubular neighborhood $N$ of $\Sig$ in
$X$ and a tubular neighborhood of $\Sig$ in $H_{g,n}$.  Using that
diffeomorphism, pull back and extend arbitrarily to $H_{g,n}$ the cone-angle
metric $\cg$, the angular 1-form $\h$,
the line bundle $\LL$ and the off-diagonal components of the adjoint connections
ad $A^0$ and ad $A^\a$. Denote these connections on $\LL$ by $B^0$
and $B^\a$ respectively. There is no obstruction to extending
$\LL$ to $H_{g,n}$ since $\S$ is a primitive class in $H_2(H_{g,n}; \Z)$.
Near $\Sig$, $B^\a = B^0 + 2i\a \h$.

Let $\cDak$ be the twisted Laplacian
acting on $\LL$-valued $k$-forms corresponding to the connection $B^\a$,
formed using the cone-angle metric $\cg$.
Consider the blow-up $\tHgn$ of $H_{g,n}$ and lift $\cDak$
to $\tHgn$. Locally near $\d \tHgn$, the 1-form $\h$ can be integrated to
(almost) an angular variable $\tth \in \R / 2\pi \Z$.  A truly angular
variable is given by $\th = \ak  \tth \,
 \in \, \R / \ak 2 \pi \Z$. We have also the distance variable $r$ and we 
choose local variables (s,t) on $\Sig$ extended to $H_{g,n}$ and lifted
to $\tHgn$. 

\begin{prop}  \label{P:Lapledge}
The operators $r^2 \cDak$  are elliptic edge for
the natural edge structure on $\d \tHgn$.
\end{prop}
\begin{proof}
This is a tedious, though straightforward, computation. Here is how to do it
for $\cDa$. Let $(x,y)$ be the Euclidean
coordinates corresponding to the polar coordinates $(r, \tth)$, i.e.
\begin{equation}   \label{E:euclpolar}
  \begin{cases}
      x = r \cos \tth  & \! = r \cos \ka \th \, , \notag  \\
      y = r \sin \tth  & \! = r \sin \ka \th.  \notag
  \end{cases}
\end{equation}
Let $d_0$ and $d_\a$ denote the covariant exterior derivatives
associated with $B^0$
and $B^\a$ respectively.
Then, in some local trivialization of $\LL$, $d_0$ takes the form
$$
 d_0 u = ( \d_s u + B_s u) \, ds + ( \d_t u + B_t u) \, dt
 +( \d_x u + B_x u) \, dx
   + ( \d_y u + B_y u) \, dy,
$$
where $B_s, \dots, B_y$ are smooth imaginary-valued functions.
Since $\cg$ is a Riemannian product near $\Sig$, one can already conclude
that the terms of $\cDa$ involving partial derivatives
in $s, t$ and $B_s, B_t$ can be grouped to form an operator which is elliptic
in $\d_s, \d_t$.  Therefore, we only need to consider the terms involving
$x, y$. Drop the $(s,t)$--terms in $d_0$ and go to $(r, \th)$--coordinates.
Then, $d_0$ becomes
$$
    d_0 u = (\d_r u +  ( B_x \cos \ka \th  +   B_y \sin \ka \th) u) \, dr
       + (\d_\th u + ( - \ka B_x r \sin \ka \th  + \ka B_y r \cos \ka \th ) u)  \, d\th \, ,
$$
which can be rewritten
\begin{equation}    \label{E:dback}
     d_0 u = (\d_r u + B_r u) \, dr  + ( \rr \d_\th u +  B_\th u ) \,  r d \th \, ,
\end{equation}
where $B_r, B_\th$ are smooth imaginary-valued functions. Hence,
$$
   d_\a u = (\d_r u + B_r u) \, dr  + ( \rr \da u +  B_\th u ) \,  r d \th .
$$
The cone-angle metric is given by $dr^2 + r^2 d\th^2$.
Therefore, the $(dr, d \th)$--components of  $\cDa$ are computed as follows:
\begin{equation}  \label{E:Laplrtha}
  \begin{split}
  - \cst d_\a \cst d_\a u  &= - \cst d_\a \{ (\d_r u + B_r u) \, r d\th  -
         ( \rr \da u +  B_\th u ) \, dr \}  \\
     &= - \{ \d_r^2 u +  ( \rr + 2B_r) \d_r u +  
       \frac{1}{r^2} \da^2 u + \frac{2}{r} B_\th \da u   \\
      & \qquad \qquad \qquad \qquad \qquad
           + ( \d_r B_r + \rr B_r + B_r^2 + \rr \da B_\th  + B_\th^2 ) u \}. 
  \end{split}
\end{equation}
From \eqref{E:Laplrtha} and the remark above regarding the $(s,t)$--terms,
it is now clear
that $r^2 \cDa$ is edge and elliptic.
\end{proof}

From the calculations in the proof of proposition \ref{P:Lapledge}, we
can carry out the analysis for $\cDa$ (Fredholmness, elliptic regularity, etc).
First, we compute the typical element
in the normal family of  $\rcDa$. Pick a point $z \in \S$ and  
coordinates $(s,t)$ centered at $z$ on $\S$, which are orthonormal at $z$.
 Lift
and extend those coordinates to $\tHgn$ in the usual way. Referring to 
equation \eqref{E:normal} in the appendix, we see
the normal operator of $\rcDa$ at $z$ is obtained by writing $\rcDa$ as
a polynomial expression in the operators $r\d_s, r\d_t, r\d_r, \d_\th$ with
smooth coefficients, setting $r=0$ in the coefficients, and viewing the
variables $(s,t)$ in the coefficients just as parameters.  Hence, at $z$,
we just set the parameters $s=t=0$.  Since the $(s,t)$--terms in $\cDa$ are
plain elliptic, after multiplying by $r^2$ and setting $r=0$ in the
coefficients, we are just going to be left with the principal 
partial $(s,t)$--symbol
of $\cDa$, multiplied by $r^2$.  As $s$ and $t$ are orthonormal at $z$, it 
implies that the $(s,t)$--terms in $N(\rcDa)_z$ are $-(r^2 \d_s^2 + r^2 \d_t^2)$.
To get the $(r, \th)$--terms, we use expression \eqref{E:Laplrtha}. Therefore,
we obtain
\begin{equation}  \label{E:normLap}
 N(\rcDa)_z = -\bigl (r^2 \d_s^2 + r^2 \d_t^2 + r^2 \d_r^2 + r \d_r +
        \da^2 \bigr) .
\end{equation}
That expression shows that all the operators in the normal family are
indeed similar.  From now on, we
won't mention the dependence on $z \in \S$ anymore.

To get the indicial operator, remark by comparing equations \eqref{E:normal}
and \eqref{E:indicial} that one just needs to drop the $(s,t)$--terms
in \eqref{E:normLap} :
\begin{equation}  \label{E:indicLap}
 \begin{split}
 I(\rcDa) &= -\bigl( r^2 \d_r^2 + r \d_r +  \da^2 \bigr) \\
          &= -\bigl( (r \d_r)^2 +  \da^2 \bigr) .
  \end{split}
\end{equation}
As in \eqref{E:indMellin}, substitute the complex variable $\z$ for $r\d_r$
in \eqref{E:indicLap}  to obtain the Mellin transform:
\begin{equation}  \label{E:MellinLap}
  I_\z(\rcDa) = -\bigl( \z^2 +  \da^2 \bigr) .
\end{equation}
The spectrum of the skew-adjoint operator $\d_\th + 2i\k$ on $S^1$,
 where the variable
$\th$ runs from $0$ to $\ak 2\pi$, is
\begin{equation}  \label{E:specS1}
 \begin{split}
   Spec \da &= \{ i(\ka m + 2\k) \, , \, m \in \Z \}  \\
              &=  \{ i\k ( \tfrac{m}{\a} + 2) \, , \, m \in \Z \}.
 \end{split}
\end{equation}
From \eqref{E:specS1}, the boundary spectrum of $\rcDa$, that is
the values of $\z$ for which  \eqref{E:MellinLap} is not invertible,
is computed to be
\begin{equation}  \label{E:specbLap}
   \sb(\rcDa) = \{ \pm \k | \ma + 2 | \, , \, m \in \Z \}.
\end{equation}
Let $\l_m = | \ma +2 |$. Since $\a \in (0, \frac{1}{2} )$,
 the set $\{ \l_m \, , \, m \in \Z \}$ is bounded
away from zero.  Indeed, the value of $\l_m$ closest to zero is given either
by $m=0$ or $m=-1$.  Let
\begin{equation}  \label{E:boundspecbLap}
   \g \: = \: \min \{ \l_m \, , \, m \in \Z \} \: = \: \min \{ 2 \, , \tfrac{1-2\a}{\a} \} .
\end{equation}

So, we have

\begin{lemma}  \label{L:specbLap}
The real part of the boundary spectrum of $\rcDa$, $\rsb(\rcDa)$, can be moved
away from zero as much as desired by taking $\k$ large enough, i.e. by taking a 
sharp enough cone-angle.  More precisely,
$$
  \rsb(\rcDa) \cap (-\k \g, \k \g) = \O .
$$
\end{lemma}

The normal operator $N(\rcDa)$ is isotropic in $(s,t)$.  To 
compute its associated Bessel-type operator $L_0$  \eqref{E:Besseltype},
substitute $i$ for $\d_s$ and
drop the $\d_t$--term:
\begin{equation}  \label{E:BesselLap}
  L_0 = -r^2 \d_r^2 - r \d_r + r^2 - \da^2
\end{equation}
The operator $L_0$ is Fredholm on $r^\de L^2 (\R^+ \times S^1)$ iff 
$\de \notin \rsb(\rcDa) + \half $ and we need to find the values of $\de$ for which
$L_0$ is invertible on $r^\de L^2$.
The kernel of $L_0$ can readily be found by separation of variables,
using \eqref{E:specS1}. 
On the $m^{th}$ eigenspace of $\da$, $L_0$ acts as
\begin{equation}  \label{E:Besseleigen}
  -r^2 \d_r^2 - r \d_r + \bigl( r^2 + \k^2 \l_m^2 \bigr) .
\end{equation}
This yields a modified Bessel equation. Its solutions are given by the
Kelvin functions (or Bessel functions of the second kind) $K_{\k \l_m}$ and
$I_{\k \l_m}$ \cite{Le}. We can 
reject the solution $I_{\k \l_m}$ which has exponential growth.
On the other hand, $K_{\k \l_m} \sim r^{-\k \l_m} $ as $r \ra 0$
 and decays exponentially as $r \ra \infty$. Therefore, 
\begin{equation}  \label{E:infdelLap}
  \ker \, L_0 = 0  \quad \text{iff} \quad \de \geq -\k \g + \tfrac{1}{2}.
\end{equation}
Next, we need to determine when coker~$L_0 = 0$. We use the fact that $L_0$ is 
self-adjoint with respect to the natural inner product in
$r^{1/2} L^2$.  The dual of $r^\de L^2$ in that inner product
is naturally identified with $r^{1-\de} L^2$.  Hence, we simply need
to know when $K_{\k \g} \notin r^{1-\de} L^2$, which gives the condition
\begin{equation}  \label{E:supdelLap}
  \text{coker} \, L_0 = 0 \quad \text{iff} \quad  \de \leq \k \g + \half.
\end{equation}

Combining lemma \ref{L:specbLap} with \eqref{E:infdelLap} and \eqref{E:supdelLap},
we have

\begin{lemma}   \label{L:Besselinv}
The Bessel-type operator $L_0$ associated to $\rcDa$ is invertible on
$r^\de L^2(\R^+ \times S^1)$ iff
$$
   - \k \g + \half \, < \, \de \, < \, \k \g + \half \, .
$$
\end{lemma}

Finally, we get to the desired result:

\begin{prop}   \label{P:invoffdag}
Let  $\a \in (0, \half)$ and let $\ell_0$ a positive integer.
Then, for $\k$ large enough, the operator
$$
  \cDa \, : \, W^p_{\ell} \ra W^p_{\ell-2}
$$
is invertible for $-\ell_0 +2 \leq \ell \leq \ell_0$.
\end{prop}

Remark that a single $\k$ won't work for all values of $\a \in (0, \half)$.
However, as will be apparent from the proof, $\k$ can be chosen so as to depend
continuously on $\a$.

\begin{proof}
By proposition \ref{P:eqedge}, the spaces $W^p_\ell$ can be identified with
certain weighted edge Sobolev spaces, the powers appearing in the weights
being bounded in terms of $\ell_0$. Lemmata \ref{L:specbLap} and \ref{L:Besselinv}
combined with theorem \ref{T:Fred} and proposition \ref{P:GLp} show that it is
possible to choose $\k$ so as to create enough elbow room in the
distribution of the indicial roots of $\rcDa$ for
$\cDa$ to be Fredholm in the given range of Sobolev differentiability.

Similarly, proposition \ref{P:wellreg} shows that elliptic regularity
between these spaces holds for $\cDa$ if $\k$ is large enough.  Since $\cDa$ is
formally self-adjoint in the cone-angle metric $\cg$,
 this implies that $\cDa$ is actually self-adjoint.
Therefore, $\cDa$ has Fredholm index $0$, and to prove the invertibility
of $\cDa$, we just need to show that $\ker \, \cDa = 0$.  
Let $u \in \ker \, \cDa$. By theorem \ref{T:Fred}, 
$u \sim r^{\k \g}$ as $r \ra 0$.  Therefore, when integrating by parts
the expression
$$
 \int_{H_{g,n}  \setminus  \S } \langle \cDa u \, , u \rangle \, = \, 0 \, ,
$$
the boundary term will converge
to $0$ as $r \ra 0$.  It follows that $u$ must be covariantly constant with
respect to the connection $B^\a$.  Since $B^\a$ is non-trivial, $u \equiv 0$.
\end{proof}

\subsection{The Laplacian for twisted connections on $\fg_E$.}  \label{SS:Laptwisted}

Recall that we have two different metrics on the manifold $\XS$, the round-metric
$g$, which extends across $\S$, and the cone-angle metric $\cg$.  Both metrics
agree away from an $\e$--neighborhood of $\S$.  Recall also
that near $\S$, the adjoint bundle $\fg_E$ splits as the sum of trivial real
bundle $i\R$ and $\LL$.  The adjoint connection ad $A^\a$ respects that splitting
near $\S$.  In fact, ad~$A^\a$ is trivial on $i\R$ and equals $B^\a$ (defined
in section \ref{SS:modeloff} ) on $\LL$.

Let $\tst$ be the Hodge star operator acting on $\fg_E$--valued forms corresponding
to the round metric for the diagonal component of $\fg_E$ and to the cone-angle
metric for the off-diagonal component near $\S$.  We can then form the formal adjoint of 
$d_{A^\a}$ in that mixed metric :
\begin{equation}
   \tda = \pm \tst d_{A^\a} \tst.
\end{equation}
The corresponding Laplacian $\tDa$ on sections of $\fg_E$ is
\begin{equation}   \label{E:Lapge}
  \tDa = \tda d_{A^\a}.
\end{equation}
In general, we will use the decoration \~{} for operators constructed using those
two different metrics for the on-- and off--diagonal components of $\fg_E$.

The following result shows that $\tDa$ behaves essentially like a usual
Laplacian between the appropriate Sobolev spaces.

\begin{prop}   \label{P:fullLapFred}
Let $\ell_0$ be a positive integer. The parameter $\k$ can be chosen large
enough so that
$$
  \tDa \, : \, \tL^p_{\ell} \ra  \tL^p_{\ell-2}
$$
is Fredholm of index $0$ for $-\ell_0 +2 \leq \ell \leq \ell_0$.
  Moreover, coker $\tDa \subset \tL^p_{\ell_0}$ and
coker $\tDa =$ ker $\tDa$.
\end{prop}

\begin{proof}
Let $X_\e$ be a small neighborhood of $X \! \setminus\!  N$ where $A^\a =A^0,
 A^0$ being the non-singular connection on $E$
chosen in section \ref{SS:fourdimsetup} (see the definition of $A^\a$ in
\eqref{E:backconn} ). Let $\D_{A^0}$ be the twisted Laplacian 
on sections of $\fg_E$ constructed using the connection $A^0$
 and the round metric $g$.
Also let $\D$ be the plain Laplacian on functions on $X$ in the round metric $g$.
  
From the construction of $\tDa$, it is clear that
\begin{equation}  \label{E:deftDa}
\tDa =
  \begin{cases}
     \D_{A^0}  &\text{over $X_\e$ ,} \\
     \cDa  &\text{on $\LL$ over $N$ ,} \\
     \D    &\text{on $i\R$ over $N$ ,}
  \end{cases}
\end{equation}
using the identification of $N$ with a neighborhood of $\S$ in $H_{g,n}$.

Let $G_{A^0}$ and $G$ be parametrices for $\D_{A^0}$ and $\D$, respectively, on $X$
and let $\G_\a$ be the inverse of $\cDa$ on $H_{g,n}$ obtained from
proposition \ref{P:invoffdag}.  A parametrix for $\tDa$ can 
be constructed in a way similar to the one described in \cite{APS} .  
Let $\phi$  be a cut-off function supported on $N$ with
$\phi \equiv 1$ on $X \setminus X_\e$, and let $\psi_1, \psi_2$
be cut-off functions supported on $N$ and $X \setminus X_\e$ respectively,
with $\psi_1 \equiv 1$ on supp $\phi$ and $\psi_2 \equiv 1$ on supp~$(1-\phi)$.
If $u$ is a section of $\fg_E$ supported on $N$, we denote by
$u_T$ and $u_D$ its off-- and on--diagonal components respectively.  

From \eqref{E:deftDa}, a parametrix for $\tDa$ is given by 
\begin{equation}   \label{E:fullparam}
 Q u =  \psi_2 G_{A^0} ((1-\phi) u) +
           \psi_1 \bigl( G((\phi u)_D) + \G_\a((\phi u)_T) \bigr) , 
\end{equation}  

Hence, $\tDa$ is Fredholm.
Moreover, $\tDa$ is formally self-adjoint for the mixed inner product
given by $g$ on $i\R$ and by $\cg$ on $\LL$.  Therefore the cokernel 
of $\tDa : \tL^p_{\ell} \ra \tL^p_{\ell-2}$ is given by the kernel of 
$\tDa : \tL^q_{-\ell+2} \ra \tL^q_{-\ell} $.  By the usual elliptic regularity,
elements of the cokernel will be smooth over $\XS$, with the
diagonal component extending smoothly across $\S$. Moreover, if $\k$ is large
enough, the off-diagonal component will belong to $W^p_2$ by proposition
\ref{P:wellreg}.  Therefore $\coker \, \tDa \subset \tL^p_{\ell_0}$.
\end{proof}

From here on,  we can consider a general connection $A = A^\a + a \in \Aa$,
where 
\begin{equation}   \label{E:genconn}
  a \, = \, \begin{pmatrix} \aD  &  a_T \\  
                       -\bar{a}_T  & -\aD 
            \end{pmatrix}
\end{equation}
near $\S$, with $\aD \in L^p_1(\Om^1(X;i\R))$ and $a_T \in W^p_1(\Om^1(\XS ; \LL))$, $p > 2$. 
On $\fg_E$--valued forms, $d_A$ acts as
\begin{equation}   \label{E:gendA} 
  d_A u = d_{A^\a} u  + 
        \begin{pmatrix}
                -2i\Im (a_T \w \bar{u}_T)   &  2(\aD \w u_T - a_T \w u_D)  \\
                2(\aD \w \bar{u}_T - \bar{a}_T) \w u_D &  2i \Im (a_T \w \bar{u}_T)
         \end{pmatrix} \, ,
\end{equation}
showing from the multiplication theorems that
 $d_A :  \tL^p_1 \ra L^p$ is a compact perturbation
of $d_{A^\a}$, see \cite{KM}.

Near $\S$, the formal adjoint of $d_A$ in the mixed metric is computed
to be 
\begin{equation}  \label{E:formadj}
  d^{\tst}_A v = d^{\tst}_{A^\a} v +
     \begin{pmatrix}
     \pm 2i \Im * ( \bar{a}_T \w \cst v_T) 
             & \pm 2i \cst ( \bar{a}_T \w  * v_D) \pm 2 \cst ( \bar{a}_{\! D} \w \cst v_T) \\
    \pm 2i \cst ( a_T \w * \bar{v}_D) \mp 2 \cst ( \aD \w \cst \bar{v}_T )  
      & \mp 2i \Im * (\bar{a}_T \w \cst v_T)
      \end{pmatrix}.
\end{equation}
The signs depend on the degree of the form $v$, but are irrelevant for our 
purpose.  Since the operator norm of $\cst$
in the round metric $g$ is uniformly bounded over $\XS$, by \eqref{E:regconeHodge},
the multiplication theorems between Sobolev spaces
show that $d^{\tst}_A : \tL^p_1(\Om^1, \fg_E) \ra L^p(\fg_E)$ is
well-defined and is a compact perturbation of $d^{\tst}_{A^\a}$.  It is at this point
that R\aa de 's techniques break down.  When writing down the adjoint of
$d_A$ in a mixed round/hyperbolic metric, one realizes that this adjoint
is not well--defined between $L^p$ spaces unless $a_T$ vanishes near the surface.

We can now
form the twisted Laplacian $\tD_A$ associated to $A$ in the mixed inner product :
$$
  \tD_A = d^{\tst}_A d_A.
$$
And we have for $\tD_A$ a result similar to proposition \ref{P:fullLapFred}:

\begin{prop}  \label{P:genconnFred}
The parameter $\k$ can be chosen large
enough so that for any singular connection $A \in \mathcal{A}^{\a,p}$, the
Laplacian
$$
  \tD_A \, : \, \tL^p_2 \ra L^p
$$
is Fredholm of index $0$.  Moreover, $\coker \, \tD_A \subset \Lpga$ and
$ \coker \, \tD_A = \ker \, \tD_A$.
\end{prop}

\begin{proof}
$\tD_A$ is a compact perturbation of $\tDa$; therefore, $\tD_A$ is also Fredholm
of index $0$.  Moreover, since $\tD_A$ is formally self-adjoint,
ker $\tD_A \subseteq$ coker $\tD_A$.  Since the dimensions have to agree,
the latter inclusion is an equality.
\end{proof}

Using proposition \ref{P:genconnFred}, we show
that ker $d^{\tst}_A$ provides a slice to the action of the gauge group
$\mathcal{G}$.

\begin{prop} \label{P:slice}
For $A \in \mathcal{A}^{\a,p}$, the image of $d_A$ is closed, and admits a closed
complement in $L^p_{1, A^\a}$, given by the kernel of $\td_A$.
\end{prop}

\begin{proof}
We essentially reproduce the proof from \cite{Morgan}. Let 
$$
 \G_A : ( \ker \, \tD_A)^{\tilde{\perp}} \cap L^p \,  \ra \,
     ( \ker \, \tD_A)^{\tilde{\perp}} \cap L^p_{2,A^\a} 
$$ 
be the inverse of $\tD_A$,
which is well-defined by proposition \ref{P:genconnFred}.  The orthogonal
complements are taken with respect to the mixed inner product.
Then the continuous maps
$$
  d_A \G_A \td_A : \tL^p_1(\Om^1, \fg_E) \ra \im \, d_A \subset \tL^p_1(\Om^1, \fg_E)
$$
and
$$
 \text{Id} - d_A \G_A \td_A :
    \tL^p_1(\Om^1, \fg_E) \ra \ker \, \td_A \subset \tL^p_1(\Om^1, \fg_E)
$$
are complementary idempotent operators, and therefore, yield a topological decomposition
$$
   \tL^p_1(\Om^1, \fg_E) = \im \, d_A \oplus \ker \, \td_A.
$$ 
\end{proof}

As in \cite{FU} proposition \ref{P:slice} can be used with the implicit 
function theorem to prove :

\begin{cor}
The space $(\mathcal{B}^{\a,p})^{*}$ 
is a Hausdorff Banach manifold.
\end{cor}

As said earlier, this was already proved by Kronheimer and Mrowka
\cite[proposition 5.6]{KM}. The main differences is that the slice
is now provided explicitly.

Proposition \ref{P:slice} also allows one to roll up the deformation complex
\eqref{E:ASDcomplex} for the ASD equation.

\begin{prop}  \label{P:ASDwrapup}
For $p$ close enough to $2$ and $A \in \mathcal{A}^{\a,p}$, the operator
\begin{equation}  \label{E:wrapcompact}
  d_A^{\tst} + d_A^+ : \tL^p_1(\Om^1(\XS ; \fg_E) \ra
               L^p(\Om^0(\XS ; \fg_E) \oplus L^p(\Om^+(\XS ; \fg_E)
\end{equation}
is Fredholm and its index is given by \eqref{E:dimASD} .
\end{prop}

\section{Laplacians on forms.}   \label{SS:Lapform}

The analysis carried out for the Laplacian on sections of $\fg_E$ in sections
\ref{SS:modeloff} and \ref{SS:Laptwisted} can be repeated for $\fg_E$--valued
forms. This will be necessary in the next chapters to study the Hodge--de~Rham
theory of singular flat connections on three-manifolds, and to develop the analysis
for ASD singular connections on manifolds with cylindrical ends.

Let $\tDak$ denote the twisted Laplacian on $\fg_E$--valued $k$--forms
formed from $A^{\a}$ over $\XS$ with the mixed metric.  
The main technical point lies in the computation of the indicial
roots of the off-diagonal component of $\tDak$.  
Recall that the off-diagonal component is studied on a model space $H_{g,n}$, from
a connection $B^\a$ pulled-back from $X$ near $\S$. The
twisted Laplacian on
$\LL$--valued $k$--forms obtained from $B^\a$ using the cone-angle metric
$\cg$ is denoted $\cDak$.  We will study $\cDaf$ in detail, because, in fact,
the indicial roots for Laplacians on forms of higher degree 
consist of the indicial roots of $\cDaf$ and $\cDa$.

Using the coordinates $(s,t,r,\th)$ on $\tHgn$ as in section \ref{SS:modeloff},
we can pick $\{ ds, dt, dr, rd\th \}$ as local generators for the pull-back of
$\Om^{*}(H_{g,n} \! \setminus \! \S)$ to $\tHgn$.  Calculations similar to proposition
\ref{P:Lapledge} show that $r^2\cDak$ is an elliptic edge differential operator.
Moreover, in the local frame $( ds, dt, dr, rd\th )$, the normal operator of 
$\rcDaf$ takes the form :
\begin{equation}    \label{E:normf}
  N(\rcDaf) \, = \,
     \begin{pmatrix}
         N(\rcDa) &0 &0  &0\\
        0 &N(\rcDa) &0  &0\\
        0 &0 &N(\rcDa)+1 &   2 \da    \\
        0 &0 & -2 \da      & N(\rcDa)+1
      \end{pmatrix}\, .
\end{equation}
This is just like the twisted Laplacian on 1-forms in 
cylindrical coordinates on the punctured 4-dimensional Euclidean space
with an axial 2-dimensional plane removed. 
Compare with \eqref{E:normLap} which gives the corresponding
Laplacian on functions.  Let $\rpDa$ denote the lower-right $2$ by $2$
matrix operator :
\begin{equation}   \label{E:normfp}
 \rpDa\,  = \,  \begin{pmatrix}
     N(\rcDa)+1 &   2 \da    \\
           -2 \da      & N(\rcDa)+1
    \end{pmatrix},
\end{equation}
in the local frame $(dr, rd\th)$.

From \eqref{E:normf}, we see that the indicial roots of $\rcDaf$ comprise those
of $\rcDa$ and those of $\rpDa$.  The indicial roots of $\rcDa$ were computed in
section \ref{SS:modeloff} and were shown to be as far away from zero as
desired by taking the parameter $\k$ large enough. The same holds for the 
indicial roots of $\rpDa$, but the proof is somewhat more involved.  
The indicial operator of $\rpDa$ is
\begin{equation}   \label{E:indfp}
  I(\rpDa) =  \begin{pmatrix}
       -(r\d_r)^2 +1 -\da^2  &2\da  \\
      -2\da & -(r\d_r)^2 +1 -\da^2
     \end{pmatrix} \, .
\end{equation}
This looks in fact like the twisted Laplacian on $1$--forms in polar coordinates on the
punctured Euclidean plane.

The Mellin transform of $I(\rpDa)$ is 
\begin{equation}   \label{Melf}
 I_\z(\rpDa) \,  = \, \begin{pmatrix}
                   -\z^2+1 -\da^2  & 2\da  \\
                -2\da & -\z^2+1-\da^2
       \end{pmatrix}.
\end{equation}
We can decompose the coefficients of $dr$ and $rd\th$ according to the spectral
decomposition of $\da$ given in \eqref{E:specS1} . The indicial roots will be
given by the values of $\z$ for which the matrix
\begin{equation}   \label{E:matreigen}
  \begin{pmatrix}
     -\z^2+1+\k^2 ( 2+ \tfrac{m_1}{\a} )  &  2i\k(2+ \tfrac{m_2}{\a})  \\
      -2i\k(2+ \tfrac{m_1}{\a})  & -\z^2+1+\k^2 ( 2+ \tfrac{m_2}{\a} )
   \end{pmatrix}
\end{equation}
is singular for some $(m_1, m_2) \in \Z \times \Z$.

Let $x=2+ \tfrac{m_1}{\a}$ and $y=2+ \tfrac{m_2}{\a}$ and remember from
\eqref{E:boundspecbLap} that $|x|, |y| \geq \g > 0$, i.e. $x,y$ are bounded
away from zero.  The determinant of the matrix \eqref{E:matreigen} is
given by
\begin{equation}   \label{E:inddet}
 \z^4 - \bigl( 2+\k^2 (x^2+y^2) \bigr) \z^2 +
    \bigl( 1 +\k^2 (x^2+y^2) + \k^4 x^2 y^2 - 4\k^2 xy \bigr) .
\end{equation}
We need to find the values of $\z$ for which \eqref{E:inddet} vanishes.  This
yields a quartic equation reducible to a quadratic equation.  The
discriminant is 
\begin{equation}   \label{E:discdet}
  \r \, = \, \k^2 \bigl( \k^2 (x^2 - y^2)^2 + 16 xy \bigr) 
\end{equation}
and the solutions are given by
\begin{equation}   \label{E:solq}
  \z^2 = \frac{ 2+\k^2 (x^2 + y^2 ) \pm \k \sqrt{\k^2 (x^2 - y^2)^2 + 16 xy}}{2}
\end{equation}

\begin{lemma}   \label{L:rspecbf}
Let
$(-\t, \t) \subset \R$. Then, for $\k$ large enough, 
$$
   \rsb (\rpDa) \cap (-\t, \t ) = \O.
$$
\end{lemma}
\begin{proof}
The proof will be divided into two cases, depending on the sign
of the discriminant $\r$.

\begin{enumerate}
 \item $\r < 0$.
     In this case, $\z^2$, and $\z$, will be complex.
  Note that $\Re (\z^2) = \dfrac{2+\k^2 (x^2 + y^2 )}{2}$
  is positive, and since $|x|, |y|$ are bounded away from zero, $\Re (\z^2)$ can be
  moved away from zero as much as desired.  The lemma now follows from the fact that
  if $z$ is a complex number with $\Re z \geq 0$, then
$$
   | \Re \sqrt{z} | \geq \sqrt{\Re z}.
$$
  \item $\r \geq 0$.  Let  $a= 2+\k^2 (x^2 + y^2 )$ and $b = \sqrt{\r}$.  There is no problem
  for the indicial roots of the form $\z = \pm \sqrt{\dfrac{a+b}{2}}$.  Let's consider the case
 where
\begin{equation}   \label{E:zab}
 \z = \pm \sqrt{\frac{a-b}{2}}.
\end{equation}
  Note that
\begin{equation}  \label{E:aqbq}
 \begin{split}
 a^2 - b^2  &= 4+4\k^2 (x^2 + y^2 ) - 16\k^2 xy + 4\k^4 x^2 y^2  \\
     &\geq C_1 \,  \k^4 (x^2 + y^2 )^2,
 \end{split}
\end{equation}
provided $\k$ is large enough, since  $|x|, |y| \geq \g$.
Moreover, direct examination shows that
\begin{equation}  \label{E:ab}
  a, b \leq C_2 \, \k^2 (x^2 + y^2 )
\end{equation}
if $\k$ is large enough.   It follows from \eqref{E:aqbq} and \eqref{E:ab} that
\begin{equation}   \label{E:finalab}
  a-b = \frac{a^2 - b^2}{a+b} \geq C_3 \, \k^2 (x^2 + y^2 ) \geq 2C_3 \, \k^2 \g^2.
\end{equation}
A glance at \eqref{E:zab} shows that  the indicial roots $\z$ can indeed be 
moved away from zero as much as desired by taking $\k$ large enough.
\end{enumerate}
\end{proof}

In view of the structure of $N(\rcDaf)$ displayed in
\eqref{E:normf}, lemmata \ref{L:specbLap} and \ref{L:rspecbf} yield

\begin{lemma}   \label{L:specbLapk}
For any interval $(-\t, \t) \subset \R$, the parameter $\k$ can be chosen
large enough so that
$$
   (-\t, \t) \cap \rsb(\rcDaf) = \O.
$$
\end{lemma}

The Bessel-type operator associated to $N(\rcDaf)$ is computed to be
\begin{equation}    \label{E:BessLapk}
  L_0 \, = \, I(\rcDaf) + r^2 \text{Id}.
\end{equation}

We need to determine when $L_0$ is invertible.

\begin{lemma}    \label{L:invBessLapk}
For any interval $(-\t, \t) \subset \R$, the parameter $\k$ can 
be chosen large enough so that $L_0$ be invertible on
$r^\de L^2$ for $\de \in (-\t +\half, \t +\half)$.
\end{lemma}
\begin{proof}
The operator $I(\rcDaf)$ is easily seen to be self-adjoint and semi-positive
with respect to the natural inner product on $r^{1/2} L^2$.  Moreover,
for $\k$ large enough, $I(\rcDaf)$ is invertible in view of lemma
\ref{L:specbLapk}.  Therefore, its spectrum is bounded away from zero.
On the other hand the multiplication operator $r^2 \text{Id}$ is
obviously positive and self-adjoint.  Therefore, Rayleigh
quotients show that the spectrum of $L_0$ is also positive
and bounded away from zero.  Hence, $L_0$ is invertible on
$r^{1/2} L^2$.  Since the kernel and cokernel of $L_0$ can jump
only for values of $\de \in \rsb + \half$, the operator $L_0$ will
stay invertible on $r^\de L^2$ for $\de$ in an interval centered around
$\half$ that can be made as large as desired by taking $\k$ large enough.
\end{proof}

\begin{prop}    \label{P:FredLap1}
Let $\a \in (0, \half)$ and let $\ell_0$ be a positive integer.  Then,
for $\k$ large enough, the operators
$$
   \cDak : W^p_\ell \ra W^p_{\ell-2}
$$
are Fredholm for $-\ell_0 +2 \leq \ell \leq \ell_0$.
\end{prop}
\begin{proof}
This follows immediately from theorem \ref{T:Fred} and lemma \ref{L:invBessLapk}.
For Laplacians on higher order froms, the only thing to remark is that the
building blocks of their normal operators will also consist of
$N(\rcDa)$ and $\rpDa$.
\end{proof}

Now, we denote  by $\tDak$ the twisted Laplacians acting on $\fg_E$-valued
k-forms, obtained using the mixed inner product.  A parametrix for $\tDak$
can be constructed in the same way as for $\tDa$ in the proof of proposition
\ref{P:fullLapFred}. Therefore, we obtain a similar result :

\begin{prop}   \label{P:fullLapFredk}
The parameter $\k$ can be chosen large
enough so that
$$
  \tDak \, : \, \tL^p_{\ell} \ra \tL^p_{\ell-2}
$$
is Fredholm of index $0$ for $-\ell_0 +2 \leq \ell \leq \ell_0$.  Moreover, 
elliptic regularity holds in that range of Sobolev
differentiability and coker $\tDak =$ ker $\tDak \subset \tL^p_{\ell}$.
\end{prop}

From the results obtained so far, it would be easy to recover and streamline
the analytical results of \cite{KM}, like regularity, compactness, etc.  For instance, 
our proof that $\im \, d_A$ is closed is much more straightforward 
\cite[lemma 5.5]{KM}. The trick
of \cite{DS} would still be needed to change the conformal structure on the 
off-diagonal component of the connection.  Instead of repeating the anlysis
of~\cite{KM} in that new framework, we are going to show how the results obtained here
can be used to develop a theory of moduli of singular connections on certain manifolds
with cylindrical ends.   The analytical techniques are inspired from those of
\cite{LM}, with a number of tricks coming from
\cite{DS}.  But  we  first need to investigate  the 3-dimensional picture.

\vfil\eject
\chapter{Three-dimensional Theory.}    \label{S:3dim}

In this chapter, we develop a Chern-Simons theory for connections
over a $3$--manifold $M$ which are singular along a knot $K$. The theory
could easily be extended to cover the more general case of a link.
The first section describes how the set-up adopted in the 
four--dimensional case can be adapted to the three-dimenional
case.   The second section defines the Chern--Simons functional $CS^\a$
as a section of a flat circle bundle on the space of connections
modulo gauge equivalence.  The critical points of $CS^\a$
are as usual flat connections and the gradient--flow lines
for $CS^\a$ can be seen as the singular ASD connections on the cylinder
$(M,K) \times \R$ studied in chpater $3$.  
The last section studies the Hessian of the Chern--Simons functional.
There, the Hodge--de Rham
theory of flat singular connections is studied, leading to a topological
condition for determining when a flat singular connection is a 
non--degenerate critical point of $CS^\a$ (proposition \ref{P:hommixed},
theorem \ref{T:homround}).

\section{Set-up.}  \label{SS:threedimsetup}
Let $K$ be an oriented knot in a homology 3-sphere $M$, endowed with some Riemannian
metric.  In this case, the set-up
described in section \ref{SS:fourdimsetup} carries over directly, but
can be somewhat simplified.  First, the knot $K$ has a natural framing
and we denote its longitude by $\l$ and its meridian by $\mu$. 
As above, we choose
a connection 1-form $\h$ for the normal circle bundle of $K$.  
We consider the trivial
$SU(2)$-bundle over $S^3$ and its trivial $U(1)$-subbundle corresponding
to the diagonal embedding of $U(1)$ in $SU(2)$.  We now  take as
the background connection $A^\a$ the  flat abelian connection
 (unique up to gauge) with 
holonomy $\exp (-2i\pi \a)$ around the meridian of $K$.  We also arrange
that the connection 1-form of $A^\a$ agrees with $\h$ over a tubular
neighborhood of $K$.  Then, for $p \geq 2$, the 
definitions \eqref{E:spaceconn} and \eqref{E:spacegauge} 
of the spaces of connections $\mathcal{A}^{\a,p} (M,K)$
and gauge
transformations $\mathcal{G}^p (M,K)$ carry over immediately, and we can form
the corresponding quotient $\mathcal{B}^{\a, p}(M,K)
= \mathcal{A}^{\a,p} (M,K) / \mathcal{G}^p (M,K)$.

The analysis carried out in the four-dimensional case carries over immediately
to the three-dimensional case.  In fact the indicial
operators associated to the twisted Laplacians on the off-diagonal
components of the connections in the cone-like metric are the
same in either dimension.  Moreover, due to the isotropy of the normal operators,
the associated Bessel-type operators are also identical.  Therefore, if we take 
a sharp enough cone-angle for the off-diagonal component of the mixed metric,
we have among other things:

\begin{prop}    \label{P:spaceconn3}
The space of irreducible singular connections $(\mathcal{B}^{\a, p})^{*}$ on $M$ is
a Hausdorff Banach manifold. A slice transverse to the action of the gauge group
$\mathcal{G}^p$  at the connection A is given by $\ker \, d_A^{\tst}$. 
\end{prop}

 We will refer to
the real algebraic variety of conjugacy classes of representations of $\pi_1(\MK)$
into $SU(2)$ as the character variety of $\MK$, and will denote it $\chi (K)$. 
The set of gauge-equivalence classes of flat connections on $\MK$ is in one-to-one
correspondence with $\chi(K)$.  By restriction to the boundary of a tubular
neighborhood of $K$, $\chi(K)$ maps into the
character variety of the torus, commonly referred to as the the ``pillowcase''.

\begin{prop}  \label{P:gaugetransfo3}
A gauge transformation $g \in \mathcal{G}^p(M,K)$ is continuous and maps $K$ into $U(1)$.
Moreover,
$$
  \pi_0(\mathcal{G}^p) \cong \Z \oplus \Z \, ,
$$
the first integer being given by the degree of $g : M \ra SU(2)$, and the second
integer by the degree of $g|_K : K \ra U(1)$.
\end{prop}

\begin{proof}
Since $W^p_2 \subset L^p_2$ and $L^p_2 \subset C^0$ by the Sobolev embedding theorem,
we see that $g$
is indeed continuous.  Now assume that for some $z \in K$, the off-diagonal
component $g_T$ of $g$ doesn't vanish.  Then $g_T$ is bounded away from
zero on a whole neighborhood of $z$ in $M$.  But this implies that
$\tfrac{1}{r^2} g_T \notin L^p$, a contradiction.

Therefore $\pi_0(\mathcal{G}) = [(M,K),(S^3,S^1)]$. Showing that the latter is 
isomorphic to $\Z \oplus \Z$ is an elementary exercise in obstruction theory, which
we present for the sake of completeness.
To show surjectivity, take any map $f : K \ra S^1$ with the desired $S^1$-degree.
Since $H^4(M,K) = 0$, $f$ can be extended to $M$ as a map into $S^3$. Moreover,
the $S^3$-degree of $f$ can be changed arbitrarily by changing $f$ on a small ball
away from $K$.  For injectivity, assume that $f$ and $g$ are two maps with
the same degrees.  It is then possible to find a homotopy
$H : [0,1] \times K \ra S^1$ between $f$ and $g$ on $K$.  The obstruction to extending
$H$ to $S^3$ lies in $H^4([0,1] \times M, \{ 0,1 \} \times M \cup [0,1] \times K)
\cong H^4([0,1] \times M, \{ 0,1 \} \times M) \cong \Z$ and is precisely given by
the difference in the $S^3$-degrees.
\end{proof}
 
\begin{cor}   \label{C:fundgroup}
$$
  \pi_1((\mathcal{B}^{\a,p})^*) \cong \Z \oplus \Z.
$$
\end{cor}

\begin{proof}
Let $\hat{\mathcal{G}}^p = \mathcal{G}^p / {\pm 1}$. Then $\pi_0(\hat{\mathcal{G}^p}) \cong \Z \oplus \Z$,
since $1$ and $-1$ are obviously homotopic.  The space $(\mathcal{A}^{\a,p})^*$
of irreducible connections  is contractible, because the space of abelian
connections is a submanifold of $\mathcal{A}^{\a,p}$ of infinite codimension
\cite{Morgan}.
Now, the action of $\mathcal{G}^p$ on
$(\mathcal{A}^{\a,p})^{*}$ (the space of irreducible connections) descends to a free
action of $\hat{\mathcal{G}^p}$ on $(\mathcal{A}^{\a,p})^{*}$.  Therefore,
$\pi_1((\mathcal{B}^{\a,p})^{*}) \cong \Z \oplus \Z$.  
\end{proof}

\section{The Chern-Simons functional.}   \label{SS:CS}
As mentioned briefly in \cite{Kr}, it is possible to define an analogue
to the Chern-Simons functional for singular connections.  However, in contrast 
to the usual Chern-Simons function, which is $\R / \Z$-valued, the functional
that we will define is in fact a section of a flat $S^1$ bundle over the
space of connections $\mathcal{B}^\a$.  The idea of defining a Chern-Simons functional
as a section of a bundle goes back to \cite{RSW} for 3-manifolds with boundary.
However, in their definition, the $S^1$--bundle is not flat.  As in the usual case,
the gradient of our Chern-Simons functional (with respect
to the flat connection on the $S^1$--bundle) will be essentially the curvature
operator.

\subsection{Definition.} 
As in \cite{DK}, the approach we are taking to define the Chern-Simons functional is by 
integrating the energy of an extension of the connection to a four--dimensional
manifold.
Let $M$ bound a four-manifold with boundary $X$ and choose a surface $\S \subset X$
so that $\S \cap M = \d \S = K$.  Extend the trivial $SU(2)$-bundle on $M$ to
a bundle $\tE$ on $X$, and the diagonal $U(1)$-subbundle to a subbundle
$\tilde{L}$ on a tubular neighborhood $N$ of $\S$ in $X$. Also, extend the form $\h$
to $N$.  Let $A$ be a singular connection 
in $\mathcal{A}^{\a,p}(M)$.  Extend $A$ to a singular connection $\tA \in \mathcal{A}^{\a,p}(X)$.
Then we define
\begin{equation}   \label{E:CSlifted}
  \wCS (A) = \frac{1}{8\pi^2} \int_{X \setminus \S } \tr \, F_{\tA} \w F_{\tA} \,
   - c_2(\tE)[X] + 2\a c_1(\tL)[\S] +\a^2 \S \cdot \S  ,
\end{equation}
where $c_1(\tL) \in H^2(N, N \cap M)$ is the relative Chern class of $\tL$ with
respect to its fixed trivialization on $N \cap M$, $c_2(\tE) \in H^4(X,M)$
is the relative Chern class of $\tE$, and the self-intersection of $\S$ is taken 
with respect to the natural framing of $K = \d \S$ in $M$.
The Chern--Weil formula \eqref{E:CW} shows that the definition \eqref{E:CSlifted}
is independent of the choices made in the extensions.  To see this,
one just needs to glue two different extensions together.  Moreover,
the same formula shows that if $g$ is a gauge transformation in $\mathcal{G}^p$, then
\begin{equation}  \label{E:gaugechange}
   \wCS (g^{*}A) = \wCS (A) + \deg \, g - 2\a \, \deg \, g|_K.
\end{equation}

Consider the trivial $S^1$-bundle  $\mathcal{A}^\a \times S^1$.  Define
a flat $S^1$-bundle over
$\mathcal{B}^\a$ by the identification
$$
     (A,z) \sim (g^{*}A, z \, \text{exp} (-4i\pi \a \, \deg \, g|_K) )  .
$$
Then, \eqref{E:gaugechange} shows that $\wCS$ descends to $\mathcal{B}^\a$
as a section of that flat bundle.  This section is our Chern-Simons functional and
will be denoted $CS^\a$.

We will give formulae for $CS^\a$, or rather $\wCS$, in terms of connection matrices,
and we  are going to compute its differential.

\begin{prop}     \label{P:CSmatrix}
Let $A \in \mathcal{A}^{\a,p}$ and let $b = \frac{1}{i} A_D - \h$, where $A_D$ is
the diagonal component of $A$.  Then
\begin{equation}    \label{E:CSmatrix}
  \wCS (A) = \frac{1}{8\pi^2} \int_{M \setminus K} \tr \, ( A \w dA +
      \tfrac{2}{3} A \w A \w A )  + \frac{\a}{2\pi} \int_K b.
\end{equation}
\end{prop}
\begin{proof}
First, note that all the terms in the RHS of \eqref{E:CSmatrix} make sense
and that it is continuous in $A$ by the Sobolev
multiplication and restriction theorems.  Therefore, we can
reason by density, and assume that $A$ is smooth on $\MK$ and that
its off-diagonal component vanishes close enough to $K$.
  
Choose $X$ and $\S$ as in the definition of $\wCS$,
but arrange that $\S$ be a disk with zero self-intersection, which can always
be done.  The extended bundles $\tE$ and $\tL$ can also be chosen to be trivial,
and let $\tA$ be an extension of $A$ to $\tE$.  For $\e$ small,
let $M_\e$ be the complement of
the $\e$-neighborhood of $K$ in $M$ and let $\S_\e$ be the boundary of 
the $\e$- tubular neighborhood of $\S$ in $X$. Then $M_\e \cup \S_\e$ is a closed
3-manifold. From the extension definition of the Chern-Simons functional for 
non-singular connections,
\begin{equation}  \label{E:limCS}
  \begin{split}
  \wCS(A)  &= \lim_{\e \ra 0} CS ( \tA |_{M_\e \cup \S_\e} )  \\
           &= \frac{1}{8\pi^2}  \lim_{\e \ra 0} ( \int_{M_\e} \tr ( A \w dA +
               \tfrac{2}{3} A \w A \w A )  +
                 \int_{\S_\e} \tr ( \tA \w d \tA + \tfrac{2}{3} \tA \w \tA \w \tA )  \\
           &= \frac{1}{8\pi^2}  \int_{M \setminus K} \tr ( A \w dA +
               \tfrac{2}{3} A \w A \w A )  +  \frac{1}{8\pi^2} \lim_{\e \ra 0}
               \int_{\S_\e} \tr ( \tA \w d \tA + \tfrac{2}{3} \tA \w \tA \w \tA ). 
 \end{split}
\end{equation}
We need to compute the last term. It is no restriction to assume that the extension
was chosen so that near $\S$,
$$
   \tA = \begin{pmatrix}  i \tb & 0 \\ 0 & -i \tb \end{pmatrix} +
             \begin{pmatrix} i\a & 0 \\ 0 & -i\a \end{pmatrix} \h \, ,
$$
where $\tb$ is an extension of $b$.  Hence,
$$
  \tr ( \tA \w d \tA + \tfrac{2}{3} \tA \w \tA \w \tA ) = -2 (\tb + \a \h) \w d \tb.
$$
Since $\tb$ is smooth,
$$
  \lim_{\e \ra 0} \int_{\S_\e} \tb \w d \tb = 0
$$
and 
\begin{equation}  \label{E:lasttermCS}
  \begin{split}
  \lim_{\e \ra 0} -2 \a \int_{\S_\e} \h \w d \tb 
               &= \lim_{\e \ra 0} -2 \a \int_{\S_\e} d ( \h \w \tb) \\
   &= 2\a \lim_{\e \ra 0} \int_{\d M_\e}  \h \w b \\
   &= 4\pi \a \int_K b .
\end{split}
\end{equation}
To go from the first to the second line, we used Stokes's theorem and
the fact that $\d \S_\e = \d M_\e$ with the 
opposite orientation.  Plugging \eqref{E:lasttermCS} into \eqref{E:limCS} 
yields the desired formula.
\end{proof}

\subsection{The differential and the critical points of $CS^\a$.}

\begin{prop}   \label{P:diffCS}
The differential of $\wCS$ at a connection $A$ is given by:
\begin{equation}    \label{E:diffCS}
  D_A \, \wCS (a) = \frac{1}{4\pi^2} \int_{M \setminus K} \tr \, (F_A \w a).
\end{equation}
\end{prop}

\begin{proof}
Reasoning by density, we can assume that near $K$, 
$$
  a =  i \begin{pmatrix} \r &0 \\ 0 &-\r \end{pmatrix} \, ,
$$
and
$$
  A =  i \begin{pmatrix} b+\a \h &0 \\ 0 &-b-\a \h \end{pmatrix},
$$
where $b, \r \in L^p_1( M \! \setminus \! K, \R )$. Then, from \eqref{E:CSmatrix},
\begin{equation}
 \begin{split}  \label{E:compdiff}
    D_A \wCS (a) &= \cpi \int_{M \setminus K}
              \tr ( a \w dA + A \w da + 2 A \w A \w a)  + \frac{\a}{2\pi} \int_K \r \\
    &= \cpi \int_{M \setminus K} \tr (2 F_A \w a - d(A \w a) )
                                               + \frac{\a}{2\pi} \int_K \r  \\
    &= \cpi \int_{M \setminus K} \tr (2 F_A \w a) -
      \cpi \lim_{\e \ra 0} \int_{M_\e} A \w a + \frac{\a}{2\pi} \int_K \r.
  \end{split}  
\end{equation}
Computing the limit as in the end of the proof of proposition \ref{P:CSmatrix}
shows that the last two terms of \eqref{E:compdiff} cancel.
\end{proof}

From proposition \ref{P:diffCS}, we see that the gradient of $\wCS$ is given by
\begin{equation}  \label{E:gradCSround}
   \nabla \wCS (A) = \frac{1}{4\pi^2} *F_A
\end{equation}
in the round metric, and by
\begin{equation}    \label{E:gradCSmixed}
   \tilde{\nabla} \wCS (A) = \frac{1}{4\pi^2}  \tst F_A
\end{equation}
in the mixed inner product.

\begin{prop}   \label{P:critCS}
The critical points of the functional $CS^\a$ are given by the gauge equivalence 
classes of flat connections on $\MK$ whose holonomy along the meridian of $K$ has trace
equal to $2 \cos 2\pi \a$.
In other words, the critical points are in one-to-one correspondence with the set of
$SU(2)$-representations $\r$ of $\pi_1(\MK)$ such that 
\begin{equation}   \label{E:meridhol}   
    \r ( [ \mu ] ) = \begin{pmatrix}  e^{-2i\pi \a} &0 \\ 0 &e^{2i\pi \a}  \end{pmatrix},
\end{equation}
modulo the adjoint action of $U(1) \subset SU(2)$.
\end{prop}
\begin{proof}
All we need to show is that any flat connection on $\MK$ 
with meridinal holonomy conjugate to \eqref{E:meridhol} can be gauge-transformed into
$\mathcal{A}^{\a,p}$.  Let $A$ be such a flat connection.

Let $T$ be a small tubular neighborhood of $K$. Since
$\pi_1(\TK) \cong \Z \oplus \Z$, it follows that $A_{T \setminus K}$ is abelian.
Therefore, after a gauge transformation
$g$ on $\TK$, $g^{*}A|_{T \setminus K}$ takes the form
\begin{equation}   \label{E:transfflat}
   g^{*}A_{T \setminus K} = \begin{pmatrix} i\a &0 \\ 0 &-i\a   \end{pmatrix} \h  +
                \begin{pmatrix} i\b &0 \\ 0 &-i\b   \end{pmatrix} dt \, ,
\end{equation}
where $t \in \R / 2\pi \Z$ is the longitudinal coordinate in  $T$, and $\b \in [0,1)$.
Since $g|_{\d T} : \d T \ra SU(2)$ is null-homotopic,  $g$ can be
extended to $M$.  From \eqref{E:transfflat}, it is clear that
the diagonal term of $g^{*}A - A^\a \in L^p_1$ and that its off-diagonal
term is in $W^p_1$.
\end{proof}

A connection matrix that takes the form given in \eqref{E:transfflat} near $\S$
will be said to be in \emph{nice form}.

We denote the set of critical points of $CS^\a$ by
$$
  \cFa (M,K) = \{ [A] \in \Bap (M,K)  \, \mid \, F_A = 0 \}.
$$
Proposition \ref{P:critCS} shows that $\cFa$ is independent of $p$.

We would now like to understand the deformation complex of the Chern-Simons
functional, to determine among other things, under which conditions a critical
point of $\wCS$ will be non-degenerate.

\section{The deformation complex.}  \label{SS:defcomp3}

From \eqref{E:gradCSround}, the Hessian of $4\pi^2 \wCS$ at a critical
point $A$ is given by $*d_A$.  Therefore, the deformation
complex for $CS^\a$ is given by :
\begin{equation}    \label{E:roundcomplexCS}
  0 \lra \tL^p_2 (\Om^0 (\MK \, ; \su))  \overset{d_A}{ \lra}
     \tL^p_1 (\Om^1 (\MK \, ; \su)) 
    \overset{*d_A}{ \lra} L^p(\Om^1 (\MK \, ; \su)) \lra 0.
\end{equation}

As usual in Floer theory, this complex is not Fredholm, the last mapping having
generally an infinite dimensional cokernel.
The complex is completed in the following way,
which has the advantage of tying in well with the anti-self-duality equation on
the cylinder :
\begin{multline}  \label{E:Floerroundcomplex}
 0 \lra \tL^p_2 (\Om^0 (\MK \, ; \su ))  \overset{d_A}{ \lra}
   \tL^p_1 (\Om^1 (\MK \, ; \su )) \oplus  \tL^p_1 (\Om^0 (\MK \, ; \su)) \\
    \xrightarrow{*d_A + d_A} L^p(\Om^1 (\MK \, ; \su)) \lra 0.
\end{multline}

\begin{prop}   \label{P:Fred3complex}
Let $A$ be a critical point of $\wCS$. Then, for $p$ larger than, but close
enough to $2$, the complex \eqref{E:Floerroundcomplex}
is  Fredholm, i.e. its homology is finite-dimensional and the image of each map
is closed.
\end{prop}

\begin{proof}
We are going to use the relation between \eqref{E:Floerroundcomplex} and the 
deformation complex of the anti-self-duality equation.
Let $(X,\S) = (M,K) \times S^1$ and endow $X$ with the product metric.  Let
$p : X \ra M$ be the natural projection and pull back $A^\a$ to $X$ with $p$.
Then, as is well-known \cite{Floer}, we have natural identifications
\begin{equation}
 \begin{align}
   \L^0(X) &\cong p^{*} \L^0(M) \, ,     \label{E:L0} \\
   \L^1(X) &\cong p^{*} \L^0(M) \oplus p^{*} \L^1(M) \, ,  \label{E:L1} \\
   \L^+(X) &\cong p^{*} \L^1(M).  \label{E:L+}
 \end{align}
\end{equation}
The identification \eqref{E:L1} is given by considering the $M$- and
$S^1$-components of a $1$-form on $X$.  The identification \eqref{E:L+} is
obtained by sending a $1$-form $\om$ to $(\om \w dt)^+$, the self-dual
part of $\om \w dt$, where $t$ is a unit-length coordinate on $S^1$.  When
the above identifications are made, the deformation complex for the ASD equation
\eqref{E:ASDcomplex} at $p^{*}A$ becomes
\begin{multline}  \label{E:ASDequiv}
 0 \lra \tL^p_2 (p^{*} \Om^0 (\MK \, ; \su ))  \xrightarrow{ \d_t + d_A}
   \tL^p_1(p^{*} \Om^1 (\MK \, ; \su )) \oplus 
     \tL^p_1(p^{*} \Om^0 (\MK \, ; \su)) \\
    \xrightarrow{ \d_t + *d_A + d_A} L^p( p^{*} \Om^1 (\MK \, ; \su)) \lra 0,
\end{multline}
which is known to be Fredholm for $p$ close enough to $2$ by the analysis
carried out in \cite{KM}.  Comparing 
the two complexes, we see that \eqref{E:Floerroundcomplex} is simply
the invariant part of \eqref{E:ASDequiv} under the action of the obvious
$S^1$ group of isometries on $X$.  Averaging over $S^1$ shows
that the homology of \eqref{E:Floerroundcomplex} is simply the $S^1$-invariant
part of the homology of \eqref{E:ASDequiv}, and is therefore finite-dimensional.
\end{proof}

Let $\mathcal{H}^i_A$ be the $i^{th}$ homology group of the complex 
\eqref{E:Floerroundcomplex}.  We are going to use an alteration
of \eqref{E:Floerroundcomplex} in the mixed metric as a prop 
to compute $\mathcal{H}^{*}_A$ :
\begin{multline}  \label{E:mixedcomplex}
 0 \lra \tL^p_2 (\Om^0 (\MK \, ; \su ))  \overset{d_A}{ \lra}
  \tL^p_1 (\Om^1 (\MK \, ; \su )) \oplus  \tL^p_1 (\Om^0 (\MK \, ; \su)) \\
    \xrightarrow{ \tst d_A + d_A} L^p(\Om^1 (\MK \,;, \su)) \lra 0.
\end{multline}
The only thing that has changed is that we are now considering the Hessian
$\tst d_A$ of $4 \pi^2 \wCS$ at $A$ with respect to the mixed inner product instead of 
the Hessian in the round metric. Let $\tcH^{*}_A$ be the homology of
the complex \eqref{E:mixedcomplex}.

\begin{lemma}  \label{L:Eulerchar}
The complexes \eqref{E:Floerroundcomplex} and \eqref{E:mixedcomplex}
have the same index.
\end{lemma}
\begin{proof}
On $X = M \times S^1$,
the deformation complex for the ASD equation \eqref{E:ASDcomplex}  and the corresponding
complex for an ASD equation in the mixed metric are homotopic, by changing
continuously the conformal structure on the off-diagonal component.  Moreover,
by results of \cite{KM}, the homotopy is through Fredholm complexes if
$p$ is close enough to $2$.  Therefore, the same is true for
\eqref{E:Floerroundcomplex} and \eqref{E:mixedcomplex} since they are
just the invariant parts of the two ASD complexes.  Hence, their indices 
are the same.
\end{proof}

\begin{lemma}  \label{L:splitH1}
For $A$ a flat connection,  $\im \, * \! d_A : \tL^p_1 (\Om^1) \ra L^p(\Om^1)$ is
orthogonal to $\im \, d_A :  \tL^p_1 (\Om^0) \ra L^p(\Om^1)$.  Therefore,
$\mathcal{H}^1_A$ splits as
$$
   \mathcal{H}^1_A = \mathcal{H}_A' \oplus \mathcal{H}_A'' \, ,
$$
where $\cH_A' \subseteq \Om^1$ and $\cH_A'' \subseteq \Om^0$.
Similarly, $\im \, \tst d_A$ is orthogonal to $\im \, d_A$ in the mixed inner product
and  $\tcH^{1}_A$ splits as
$$
    \tcH^1_A = \tcH_A ' \oplus \tcH_A '' .
$$
\end{lemma}

Remark that in fact, $\mathcal{H}_A'$ is the first homology group of the complex
\eqref{E:roundcomplexCS} and, as such, can be seen as the kernel of the
Hessian of $CS^\a$.

\begin{proof}
Since $A$ is flat, integration by parts shows that 
$$
  \int_{M \setminus K} \langle d_A u , *d_A v \rangle = 0
$$
if the off-diagonal components of $u \in \Om^0( \MK; \fg_E)$ and
$v \in \Om^1( \MK; \fg_E)$  vanish near $\S$.
The result follows then by density.
The same proof applies in the mixed inner product.
\end{proof}

\begin{prop}  \label{P:hommixed}
The homology groups of \eqref{E:mixedcomplex} are given by
\begin{equation}
\begin{align}
   \tcH^0_A  & \cong H^0(\MK ; \ad \, A) ,  \notag \\
   \tcH^1_A  & = \tcH'_A  \oplus \tcH''_A  \notag \\
    & \cong     \ker \, (H^1 ( \MK; \ad \, A) \ra H^1( \mu; \ad \, A))
           \oplus  H^0(\MK; \ad \, A),  \notag \\
   \tcH^2_A  & \cong \tcH'_A  .  \notag
\end{align}
\end{equation}
Therefore, the index of the complex \eqref{E:mixedcomplex}
is zero.
\end{prop}
\begin{proof}
First, assume that $A$ has been gauge-transformed, so that near $\S$, $A$
is in nice form, as given in \eqref{E:transfflat}.  Then, the analysis carried out
for the Laplacians associated to the background connection $A^\a$ applies
\emph{verbatim} to $A$, and so proposition \ref{P:fullLapFredk} applies.

Now, $\ker \, d_A = \tcH^0_A = 0$ if $A$ is irreducible.  If $A$ is abelian
and $u \in \ker \, d_A$, then $u$ is smooth and 
$u = i \, \text{diag} \, ( a, -a )$ near $\S$,
for some $a \in \R$. This shows that $u \in\tL^p_2$, and thus yields
the isomorphisms for $\tcH^0_A$ and $\tcH''_A$.
 
Lemma \ref{L:splitH1} shows that $\tcH_A' = \ker \, d_A / \im \, d_A$ and therefore,
the slice theorem (proposition \ref{P:spaceconn3} ) 
implies that $\tcH_A' \cong \ker \, d_A \cap \ker \td_A$. It follows
that $\tcH_A' \subseteq \ker \, \tD_{A,1}$.  Moreover, if the cone-angle is sharp 
enough, integration by parts shows that $\ker \, \tD_{A,1} \subseteq \tcH_A'$. Now,
$H_1( \mu; \ad \, A) \cong \R$, since $A$ is abelian but non-central on $\mu$,
and a generator for $H_1( \mu; \ad \, A)$ is given by the diagonal lifting of
$\mu$.   Since the diagonal component $\om_D$ of an element
$\om \in \ker \, \tD_{A,1}$ is 
in  $C^\infty (M)$ by elliptic regularity, we see that
$$
  \int_\mu  \om_D = 0
$$
by shrinking the meridian $\mu$ towards $K$. Therefore, $\ker \, \tD_{A,1}
\subseteq  \ker (H^1 ( \MK; \ad \, A) \ra H^1( \mu; \ad \, A))$.

Consider the natural map $\phi : \tcH_A' \cong \ker \, \tD_{A,1} \lra
\ker \, (H^1 ( \MK ; \ad \, A) \ra H^1( \mu; \ad \, A))$.  We need to show that 
$\phi$ is an isomorphism.  First, we show that $\phi$ is surjective by a
Mayer--Vietoris argument.  Let $\om$ be a smooth representative of some class in
$\ker  (H^1 ( \MK; \ad \, A) \ra H^1( \mu; \ad \, A))$.  If $N$ is a tubular
neighborhood of $K$, then $\ker \, (H^1 ( \NK ; \ad \, A) \ra H^1( \mu; \ad \, A))
\cong \R$ and is generated by the $\su$-valued 1-form
$\l =i$ diag $( dt, -dt) \in \tL^p_{1}$,
where $t$ is a coordinate on $K$, extended to $N$.
The form  $\om$ restricted to $\NK$ is therefore cohomologous to
$a\l$ for some $a \in \R$, i.e. $\om - a \l = d \xi$  for some 
section $\xi$ on $\NK$.  Hence,
$\bar{\om} = \psi a \l + (1-\psi) (\om - d \xi)$, where $\psi$ is a cut-off
function supported on $N$,  is cohomologous to $\om$ and
belongs to $\tL^p_{1}$.

To show that $\phi$ is injective, we adapt the technique
of \cite{APS}. Let $\pi : \tM \ra M$ be the oriented real
blow-up of $M$ along $K$ and pull back $A$ to $\tM$.
Assume that $\om \in \ker \, \tD_{A,1}$ is exact.
We have $\pi^{*}\om \in L^p_1(\tM)$, as is  easily checked in local coordinates
and $\om = d_A \s$ for some form $\s \in L^p_2(\tM) \subseteq r^{-1/q} L^p_2(M)$,
by Hodge theory on manifolds with boundary. Then, using the mixed metric to define
$L^2$-norms on $M$,
\begin{equation}  \label{E:injH}
\begin{split}
  \| \om \|^2_2 &= \lim_{\e \ra 0} \int_{M_\e} \tr (\tst \om \w \om)  \\
        &= \lim_{\e \ra 0} \int_{M_\e} \tr ( \tst \om \w d_A \s - d_A \tst \om \w \s) \\
        &= \lim_{\e \ra 0} \int_{M_\e} \tr ( d_A (\tst \om \w \s ) ) \\
        &= \lim_{\e \ra 0} \int_{\d M_\e} \tr (\tst \om \w \s) \\
        &= 0 \, ,
\end{split}
\end{equation}
given the asymptotics of $\om$.  To go from the first to the second line of
\eqref{E:injH}, notice that the last term of the second line 
is zero since $d_A^{\tst} \om =0$.
Hence $\om \equiv 0$ and we have proved the injectivity of $\phi$.

For the last homology group $\tcH^2_A$, remark by duality that
\begin{equation}   \label{E:H2dual}
 \tcH^2_A \cong \coker ( \tst d_A + d_A ) 
       \cong  \{ \om \in L^q \, \mid \, d_A^{\tst} \om
                     = \tst d_A \om = 0 \}.
\end{equation}
By elliptic regularity, this implies that $\coker ( \tst d_A + d_A ) =
\ker \, \tD_{A,1} \subset L^p$. Hence, $\tcH^2_A \cong \tcH_A'$.

When the connection matrix of $A$ is not in nice form,
$\tcH^0_A = \ker \, d_A$ and 
$\tcH^1_A = \ker \, d_A / \im \, d_A \oplus \ker \, d_A$ can be computed
from above by gauge equivariance.  Formula \eqref{E:H2dual} still holds
by duality and so does elliptic regularity for $\tD_{A,1}$.  Therefore, the above
reasoning goes through.
\end{proof}

The same identifications hold for the deformation complex
in the round metric \eqref{E:Floerroundcomplex}.

\begin{thm}  \label{T:homround}
For $A$ a flat connection,
$$
  \cH^{*}_A \cong \tcH^{*}_A .
$$
\end{thm}
\begin{proof}
First, notice that $\cH^0_A = \ker \, d_A = \tcH^0_A$ and
$\cH^1_A = \ker \, d_A / \im \, d_A \oplus \ker \, d_A = \tcH^1_A$. 
To compute $\cH^2_A$, we are using the fact that the complex 
\eqref{E:Floerroundcomplex} is gauge-equivariant.  Hence, we can 
assume that $A$ is in nice form.

By duality,
$$
  \cH^2_A \cong \coker ( * d_A + d_A ) 
       \cong  \{ \om \in L^q \, \mid \, d_A^{*} \om
                     = * d_A \om = 0 \}. 
$$
The group $\cH^2_A$ can be identified with the $S^1$-invariant part
of $\coker \, d_A^+$ on $( \MK ) \times S^1$. Now, the proof of \cite[lemma 5.4]{KM}
shows that $\ker (d_A^{*} + *d_A) \cap L^q \subset L^p$, i.e. a limited
version of elliptic regularity holds.  Hence, we can assume that $p=q=2$.
Let $\om \in \ker (d_A^{*} + *d_A) = \cH^2_A$.  By the standard elliptic regularity for
the diagonal component $\om_D$ of $\om$, we see that $\om_D \in C^\infty(M)$, and thus
that $[ \om ] \in \ker  (H^1 ( \MK; \ad \, A) \ra H^1( \mu; \ad \, A))$ by
reasoning as above. In view of lemma \ref{L:Eulerchar}, proposition \ref{P:hommixed} 
and the identifications already obtained, it is enough to show that any class in
$\ker  (H^1 ( \MK; \ad \, A) \ra H^1( \mu; \ad \, A))$ has an $L^2$ harmonic
representative.  The Mayer--Vietoris argument presented in the proof of 
proposition \ref{P:hommixed} shows that any such class has an $L^2$ representative.
The fact that the representative can then be chosen to be harmonic is the
content of a theorem of de Rham-Kodaira \cite{dR}.
\end{proof}

Theorem \ref{T:homround} and proposition \ref{P:hommixed} are very important because
they give topological identifications for the homology of the deformation
complexes.  In fact $\cH_A' \cong \tcH_A'$ can be seen as the (Zariski)
tangent space of $\cFa$ at $A$.  The above results express the fact that
this tangent space is just the restriction of the tangent space of the 
character variety $\chi(K)$ of $\MK$ at $A$ (identified with
$H^1(\MK ; \ad \, A)$ ) to the space of connections whose meridinal trace is
fixed by $\a$. This is what should intuitively have been expected. 

In a way similar to \cite{DFK}, we adopt the following conventions :

A flat connection $A \in \mathcal{A}^{\a,p}$ is said to be \emph{acyclic} if 
  the homology of either deformation complex vanishes.  It is said to
be \emph{isolated} if $\cH'_A = \tcH'_A = 0$.

The following is a well-known fact \cite{K} :

\begin{prop}  \label{P:abisolated}
The abelian connection with meridinal holonomy given by \eqref{E:meridhol} is isolated
iff $e^{4 i \pi \a}$ is not a root of the Alexander polynomial of $K$.
\end{prop}

In general, a non-abelian flat connection will be isolated, and thus acyclic,
if it lies on
a smooth arc of the character variety of $K$, which is transverse to the
line fixing the meridinal holonomy to be as in \eqref{E:meridhol} in the
pillowcase.

When flat connections are not isolated, the usual technique, found for
instance in \cite{TaubesCasson} or \cite{Herald}, consists of 
altering the Chern--Simons functional by so--called ``holonomy
perturbations''.  For a generic such perturbation, the critical points
of the perturbed Chern--Simons functional are isolated.  Theorem
\ref{T:homround} shows that the technique of holonomy perturbations
would also work to isolate singular flat connections.  The proof
presented in \cite{TaubesCasson} applies almost without change.
Moreover, it can be assumed that the loops supporting the
holonomy perturbations are homologically trivial in the knot
complement $\MK$.  The picture  emerging for the moduli space of perturbed
flat connections is completely similar to the
one studied in \cite{Herald} for manifolds with boundary.  The manifold
with boundary to consider here is of course the homology sphere $M$ with
a small open tubular neighborhood of the knot $K$ removed.

\vfil\eject
\chapter{Singular connections on manifolds with cylindrical ends.}

This chapter studies moduli spaces of singular ASD connections on
manifolds with cylindrical ends.  These can be seen as the spaces 
of flow lines for the gradient flow of the Chern--Simons functional
and they determine the boundary operator in Floer homology.  We are not
going as far as defining a Floer homology in this thesis, but this chapter
can be seen as a step in that direction.  The first section establishes
the basic definitions.  As usual in Floer homology, one needs to work
with Sobolev spaces which are exponentially weighted in the cylindrical
direction.  This is in fact necessary to the analytical treatment of
the case of the abelian connection.  The second section studies the basic
operators in the translation--invariant case, which serves as a local model
for the analytical theory on the ends.  The techniques used are reminiscent
of \cite{LM}.  The operators are first studied in the mixed metric and then,
the conformal structure for the off--diagonal component is changed
to the round metric thanks to the trick from \cite{DS}.
The third section presents the global theory.  One then gets
smooth moduli spaces of singular ASD connections, whose dimensions
are given by the Fredholm indices of some differential operators.
The last section says a few words on how those indices change under
gluing or when a gauge transformation is applied.  Those formulae
show that, at least as far as dimensions are concerned, the picture
is similar to the usual Floer homology.

\section{Set-up.}  \label{SS:setupcyl}

\subsection{The background connection.}
Let $W$ be a Riemannian manifold with $n$ cylindrical ends modeled on $M_j \times \R^+$,
$j = 1, \dots, n$,
where $M_j$ is a homology $3$-sphere.  Assume that the metric on $W$ is a product
metric on each end. Let
$\S \subset W$ be a surface, also with cylindrical ends, where we assume that the
ends of $\S$ are given by
$\S \cap (M_j \times \R^+) = K_j \times \R^+$, where $K_j$ is a knot in $M_j$.
Let $N$ be a tubular neighborhhod of $\S$ such that $N \cap (\MjR) = N_j \times
\R^+$, where $N_j$ is a tubular neighborhood of $K_j$ in $M_j$.  Let $\a \in (0,
\half)$. For the following, we will assume that for every knot $K_j$, the flat
connections in
$\cFa (M_j, K_j)$ are isolated.  This means that $e^{4 \pi i \a}$ is not a root of the
Alexander polynomial of $K_j$ and that the non-abelian flat connections are
acyclic.  The torus knots are examples of knots where the flat connections
in $\mathcal{F}^\a$ are isolated when $e^{4\pi i \a}$ is not a root of the
Alexander polynomial \cite{K}.

Let $E$ be an $SU(2)$--bundle on $W$, and fix trivializations of
$E|_{M_j \times \R^+} \cong  (M_j \times \R^+) \times \C^2$ on the ends
of $W$. Similarly, pick a splitting $E = L \oplus L^{-1}$ on N and fix
trivializations $L|_{N_j \times \R^+} \cong  (N_j \times \R^+) \times \C$. Consider
(isolated) flat connections $A_j \in \mathcal{F}^\a (M_j,K_j)$  in nice form for
$(M_j, K_j)$. Choose now an angular $1$-form
$\h$ around $\S$ on $\WS$ which agrees with the connection matrix
of $A_j$ on $\MKj$. 
Pick a background connection $A^\a$ on $E$ as in \eqref{E:backconn}, arranging on
each end that $A^\a|_{M_j \times \R^+} = A_j$.  Following the terminology
of \cite{DFK}, we call such a pair of bundles $(E,L) \ra (W, \S)$ with a fixed
background connection $A^\a$ which is flat on the ends an \emph{adapted bundle pair}. 

\subsection{Weighted spaces.}
Let $t_j \in (0, \infty)$ be the $\R^+$--coordinate on $M_j
\times
\R^+$ and  let $w_\de : W \ra \R, \de = (\de_1, \dots , \de_n) \in \R^n$, be smooth
positive functions such that
$$
  w_\de (x, t_j) = \exp(\de_j t_j), \qquad (x,t_j) \in \MjR.
$$
Define $\tL^p_{k,\de} (\Om^m(\WS ; \fg_E))$ to be the space of distributional
sections  $u$ of $\L^m(\WS ; \fg_E)$ such that $w_\de u$ is of class $L^p_k$ and the
off-diagonal component of $w_\de u$ near $\S$ is in $W^p_k( N \! \setminus \! \S)$.  We may
drop the subscript
$\de$ when $\de = 0$.

In the following definitions, we assume that $\de_j \geq 0$, with strict inequality
if $A_j$ is abelian.
We define the space of singular connections in the usual way :
\begin{equation}
  \mathcal{A}^{\a,p}_{\de}(W,\S) (A_1, \dots, A_n) = A^\a +
    \tL^p_{1,\de} (\Om^1( \WS ; \fg_E)). 
\end{equation}
The \emph{small gauge group} is
\begin{equation}
  \tcG^p_{\de}(W, \S) 
     = \{ g \in \text{Aut} (E) \, \mid \,
           g - 1 \in \tL^p_{2,\de} ( \WS , \text{End} \, E) \} ,
\end{equation}
and we have the corresponding quotient
\begin{equation}
 \tcB^{\a,p}_{\de}(W, \S)(A_1, \dots , A_n) = \mathcal{A}^{\a,p}_{\de} / \tcG^p_{\de}.
\end{equation}
The corresponding space of anti-self-dual connections is denoted
\begin{equation}
  \tM^{\a,p}_\de (W, \S) (A_1, \dots , A_n) = \{ \, [A] \in \tcB^{\a,p}_\de (W, \S)
     \, \mid  \,       F_A^+ = 0    \, \}.
\end{equation}
The Lie algebra of $\tcG^p_\de$ is $\tL^p_{2,\de} (\Om^2( \WS ; \fg_E ))$.
As usual, a $*$-superscript will mean that we are talking only about the subspaces
of irreducible connections.

In fact, $\tcG^p_\de$ is not the whole gauge group acting on $\mathcal{A}^{\a,p}_\de$. Instead,
as usual in Floer theory \cite{Floer} \cite{DFK}, one should consider the group
\begin{equation}
  \mathcal{G}^p_\de (W, \S)(A_1, \dots , A_n) = \{ g \in \text{Aut} (E) \, \mid \,
           \nabla_{A^\a} g \in \tL^p_{1,\de} ( \WS , \text{End} \, E) \}
\end{equation}
In fact, $\tcG^p_\de$ is a normal subgroup of $\mathcal{G}^p_\de$.  Let $\G_i$ be the isotropy
subgroup of the flat connection $A_i$, that is $\G_i = \{ \pm 1 \}$ is $A_i$ is irreducible
and $\G_i = S^1$ if $A_i$ is abelian. Then, 
\begin{equation}  \label{E:gquotient}
  \mathcal{G}^p_\de / \tcG^p_\de \cong \G_1 \times \dots \times \G_n.
\end{equation}
We can form 
\begin{equation} \label{E:goodconn}
 \begin{split}
  \Bap_\de (W, \S)(A_1, \dots, A_n) &= \mathcal{A}^{\a,p}_\de /  \mathcal{G}^p_\de \\
           &=\tcB^{\a,p}_\de / \, \G_1 \times \dots \times \G_n
 \end{split}
\end{equation}
and
\begin{equation}  \label{E:gASD}
 \begin{split}
   M^{\a,p}_\de (W, \S) (A_1, \dots , A_n) &= \{ \, [A] \in \Bap_\de (W, \S ) \, \mid \,
                  F_A^+ = 0    \, \} \\
              &=\tM^{\a,p}_\de / \, \G_1 \times \dots \times \G_n \, .
 \end{split}
\end{equation}

\section{The translation-invariant case.}
Let $(W, \S) = (M, K) \times \R$ , $\pi : W \ra M$ be the
projection and let $t$
be the axial coordinate. Consider an isolated flat connection 
$A \in \mathcal{F}^{\a,p}(M,K)$ in nice form and pull it back to $(W, \S)$.  This is
going to serve as our model for the ends. In what follows,
we assume that a sharp enough cone-angle has been chosen for
the off-diagonal component of the mixed metric, so that proposition 
\ref{P:fullLapFredk} holds for a large range of Sobolev differentiability (determined by
the integer $\ell_0$ in proposition \ref{P:fullLapFredk}).
 
\subsection{The Laplacians.}
Let $\tDM_{A,k}$ and $\tDW_{A,k}$ be the
twisted Laplacians on 
$k$-forms formed with the mixed inner product on $M$ and $W$ respectively.
Considering separately the $M$-- and $\R$-- components, the bundles of forms on W
split naturally as
\begin{equation}  \label {E:splitform}
   \L^k (W) \cong \pi^{*} \L^k(M) \oplus \pi^{*} \L^{k-1}(M).
\end{equation}
In those splittings,
\begin{equation}  
\begin{align}
  \tDW_{A,0} &= \tDM_{A,0} - \d_t^2 ,  \notag \\
  \tDW_{A,k} &= \begin{pmatrix}   \tDM_{A,k}-\d_t^2  &0 \\
                                     0 &\tDM_{A,k-1}-\d_t^2  \end{pmatrix},
      \qquad 1 \leq k \leq 3,   \label{E:Lapcylinv} \\
  \tDW_{A,4} &= \tDM_{\a,3} - \d_t^2,  \notag
\end{align}
\end{equation}
where we have identified $\tDM_{A,k}$ with its lifting to $W$.  When $A$
is acyclic, the operators $\tDM_{A,k}$ are invertible on $L^2$ and their
spectra are positive.  Strictly speaking, acyclicity means only the invertibility
of $\tDM_{A,0}$ and $\tDM_{A,1}$, but $\tDM_{A,k} = \pm \tst \tDM_{A, 3-k} \tst$.

\begin{prop}   \label{P:invLapcyl}
If $A$ is acyclic, the Laplacians $\tDW_{A,k} : \tL^p_\ell \ra \tL^p_{\ell-2}$ are
invertible for $-\ell_0 +2 \leq \ell \leq \ell_0$.
\end{prop}

\begin{proof}
The technique used here is similar to the one in \cite{LM}.  We will give the proof
only for $\tDW_{A,0}$, the other cases being entirely similar.
Conjugate $\tDW_{A,0}$ by the Fourier transform in the axial direction.
If $\t$ is the variable dual to $t$, one gets
\begin{equation}  \label{E:LapFT}
    \hD_0 (\t) = \tDM_{A,0} + \t^2 ,  
\end{equation}
The operator $\hD_0(\t)$ is invertible on $\tL^2(\MK)$ for all $\t \in \R$, and the
operator norm of $L(\t) = (\hD_0(\t))^{-1} : \tL^2 \ra \tL^2$ is uniformly bounded in
$\t$ (by the inverse of the first eigenvalue of $\tDM_{A,0}$).
Therefore, $\tDW_{A,0} : \tL^2 \ra \tL^2$ is invertible, with inverse
given by
\begin{equation}   \label{E:FinvLap}
 (G \, f)(t) = (\tDW_{\a,k})^{-1} (f) (t) = \frac{1}{\sqrt{2\pi}}
          \int_\R  L(\t) \hat{f} (\t) e^{it\t} \, d\t 
\end{equation}  
The family $L(\t) = (\hD_0(\t))^{-1}, \t \in \R,$ is an operator-valued
Fourier multiplier.  By a theorem of Mikhlin--Stein  \cite{Ty}, it maps 
$L^p(\R, \tL^2(M))$ into itself provided
\begin{equation}   \label{E:Mikhlin}
  \| \d_\t L(\t) \| \leq C \, ( 1+ | \t | )^{-1}.
\end{equation}
But $\d_\t L(\t) =  -2 \t L(\t)^2$, and \eqref{E:Mikhlin} is therefore
obviously satisfied.
Let $p \geq 2$.  From what precedes, if $f \in \tL^p(\WS) \subseteq L^p(\R,\tL^2(M))$
and $u = Gf$, then $u \in L^p(\R, \tL^2(\MK))$.

Consider the strips $B_n = [n-1,n+1] \times M$ and $B_n^+ = (n-2,n+2) \times M,
n \in \Z$.  By interior elliptic regularity (proposition \ref{P:intellestim}),
\begin{equation}  \label{E:strip}
 \begin{split}
  \| u|_{B_n} \|^p_{\tL^p_2} & \leq  C \bigl( \| f|_{B_n^+} \|^p_{\tL^p} 
          + \| u|_{B_n^+} \|^p_{\tL^2}   \bigr)  \\
    & \leq  C \bigl( \| f|_{B_n^+} \|^p_{\tL^p} 
          + \| u|_{B_n^+} \|^p_{L^p((n-2,n+2), \tL^2(\MK))} \bigr),
 \end{split}
\end{equation}
since $L^p((n-2,n+2), \tL^2(\MK)) \subset \tL^2(B_n^+)$. Summing \eqref{E:strip}
over all the integers $n \in \Z$  yields
\begin{equation} 
\begin{split}
  \| u \|^p_{ \tL^p_2(\WS)} & \leq C \bigl( \| f \|^p_{\tL^p(\WS)} +
 \| u \|^p_{L^p(\R, \tL^2(\MK))}    \bigr)  \\
     & \leq C \bigl( \| f \|^p_{\tL^p(\WS)} + \| f \|^p_{L^p(\R, \tL^2(\MK))} \bigr) \\
     & \leq C \| f \|^p_{\tL^p(\WS)}.
\end{split}
\end{equation}
Therefore, $\tDW_{A,0} : \tL^p_2 \ra \tL^p$ is an isomorphism for $p \geq 2$.
For $p < 2$ and other orders of differentiability, use elliptic regularity,
duality and interpolation.
\end{proof}

Let $\tL^p_{k,\de} (W) = e^{-t\de} \tL^p_k (W)$.  The adjoint of $d_A$ in the weighted
mixed inner product of $\tL^2_\de$ is given by
\begin{equation}   \label{E:wadj}
  d_A^{\tst_{\de}} = e^{-2\de t} d_A^{\tst} e^{2\de t}.
\end{equation}
We also need a slice theorem in the weighted case. As usual, it will be
obtained from the Laplacian. We now investigate the weighted Laplacians,
given by
\begin{equation}   \label{E:wLap}
 \tDW_{A,k,\de}\, = \, d_A^{\tst_{\de}} d_A + d_A d_A^{\tst_{\de}}:
      \tL^p_{\ell,\de} (\Om^k( \WS; \fg_E))
      \ra \tL^p_{\ell-2,\de} (\Om^k( \WS; \fg_E)).
\end{equation}

The operator $\tDW_{A,k,\de}$ is similar to the following operator
$L_{A,k,\de}$ acting on unweighted spaces:
\begin{equation}  \label{E:opL}
  L_{A,k,\de} \: = \: e^{\de t} \, \tDW_{A,k,\de} \, e^{-\de t} 
     \: : \tL^p_{\ell} \ra \tL^p_{\ell-2}.
\end{equation}
It is easily computed that
\begin{equation}   \label{E:equivL}
\begin{split}
   L_{A,k,\de} &= \tDW_{A,k} + \de^2  \\
                     &=
    \begin{cases} ( \tDM_{A,0} + \de^2 ) - \d_t^2 \qquad \qquad 
                        \quad      \qquad \qquad \text{if} \, k=0, \\
     \begin{pmatrix}   \tDM_{A,k} + \de^2   &0 \\
                                     0 &\tDM_{A,k-1}+ \de^2 \end{pmatrix}
       -\d_t^2 \qquad \text{if} \,  1 \leq k \leq 3,    \\
                  ( \tDM_{A,3} + \de^2 ) - \d_t^2  \qquad \qquad 
                        \quad \qquad \qquad \text{if} \,  k=4.
    \end{cases}
\end{split}
\end{equation}
Notice that the normal operator for the off-diagonal component of $r^2 (\tDM_{A,k} +
\de^2)$ is the same as for $r^2 \tDM_{A,k}$. Therefore, all the previous analysis
applies to $\tDM_{A,k} + \de^2$, including elliptic regularity.

\begin{prop}  \label{P:invwLap}
If $A$ is irreducible, the operators \eqref{E:wLap} are invertible for all $\de$.
If $A$ is abelian, the operators \eqref{E:wLap} are invertible
for all $\de \neq 0$, and also for $\de = 0$ if $k=2$.
\end{prop}
\begin{proof}
We simply need to show that \eqref{E:opL} is invertible.  To do this proceeds
as in the proof of proposition \ref{P:invLapcyl}, using \eqref{E:equivL}.
Note that the first terms of the RHS of \eqref{E:equivL} are always positive, since the
operators
$\tDM_{A,k}$ are positive if $A$ is irreducible or $k=2$, and $\de \neq 0$ if $A$ is
abelian and $k \neq 2$.
\end{proof}

\subsection{The ASD operator.} 
We will follow closely the technique of changing the conformal structure from \cite{DS} to
study the ASD operator in the round metric.  
We first assume that the connection $A$ is irreducible.
The case of an abelian connection
will have to be treated separately.

Let
$\Om^{\tp}(\fg_E)$ be the space of
$\fg_E$--valued
$2$--forms whose diagonal component is self-dual in the round metric and whose
off-diagonal component is self-dual in the cone-angle metric.   Define the 
operator 
$$
   d_A^{\tp} = \tfrac{1+\tst}{2} d_A : \Om^1(\fg_E) \ra \Om^{\tp}(\fg_E)
$$
as the linearization of the ASD operator in the mixed metric. Similarly, define
$\Om^{\tm}$ and $d_A^{\tm}$.  
Proposition \ref{P:invLapcyl} implies that the operator 
\begin{equation}      \label{E:wrapmixed}
  d_A^{\tst} + d_A^{\tp} : \tL^p_1(\Om^1(\WS; \fg_E)) \ra \tL^p(\Om^0(\WS ; \fg_E)) \oplus
            \tL^p(\Om^{\tp}(\WS ; \fg_E))
\end{equation}
is invertible.  Indeed, composing on the left and on the right with 
$d_A + d^{\tst}_A$ yields the invertible Laplacians. The operator
\eqref{E:wrapmixed} comes from rolling up the deformation complex for 
the ASD equation in the mixed metric.
But what we would like to show is that the operator
\begin{equation}   \label{E:wrapround}
 D_A = d_A^{\tst} + d_A^+ : \tL^p_1(\Om^1(\WS; \fg_E)) \ra \tL^p(\Om^0(\WS ; \fg_E)) \oplus
            \tL^p(\Om^+(\WS ; \fg_E))
\end{equation}
is invertible (we are now working with the rolled-up version of the 
ASD deformation complex
in the round metric for all the components of $\fg_E$).

First, let $T_A = \ker \, d_A^{\tst} : \tL^p_1(\Om^1) \ra \tL^p$.
Then, $d_A^{\tp} : T_A \ra \tL^p (\Om^{\tp})$ is an isomorphism. Let
$Q = (d_A^{\tp})^{-1} :  \tL^p (\Om^{\tp}) \ra T_A$ and let $S = d_A^{\tm} Q$.

\begin{lemma}  \label{L:isometry}
The operator $S : \tL^2(\Om^{\tp}) \ra \tL^2(\Om^{\tm})$ is an isometry.
\end{lemma}
\begin{proof}
By density, it is enough to show that for $u \in C^{\infty}_0 (\Om^1(W; \fg_E))$,
with the off--diagonal component of $u$ supported away from $\S$, we have 
$\| d_A^{\tp} u \|_{\tL^2} = \| d_A^{\tm} u \|_{\tL^2}$. Now,
\begin{equation}
 \begin{split}
   \| d_A^{\tp} u \| ^2_{\tL^2} - \| d_A^{\tm} u \|^2_{\tL^2}
    &= \int_{W \setminus \S} \tr \, (d_A^{\tp} u \w \tst d_A^{\tp} u 
          - d_A^{\tm} u \w \tst d_A^{\tm} u ) \\
   &= \int_{W \setminus \S} \tr \, (d_A u \w d_A u ) \\
   &= 0 ,
 \end{split}
\end{equation}
by integration by parts, since $A$ is flat.
\end{proof}

\begin{lemma}  \label{L:isomLp}
The operator norm of $S : \tL^p(\Om^{\tp}) \ra \tL^p(\Om^{\tm})$ tends to $1$ as
$p \ra 2$.
\end{lemma}
\begin{proof}
Combine lemma \ref{L:isometry} with the Marzinkiewicz interpolation theorem
\cite{Ty}.
\end{proof}

Recall from section \ref{SS:cone} that the change of conformal structure from
the cone-angle to the round metric is encoded in a map
$\mu : \Om^{\tm} \ra \Om^{\tp}$.  Let $\tmu: \Om^{\tm}(\fg_E) \ra \Om^{\tp}(\fg_E)$
be the map which is zero for the diagonal component of the $\fg_E$ and equal
to $\mu$ on the off-diagonal component of $\fg_E$.  Then, 
the operator $d_A^+ : T_A \ra \Om^+$ is equivalent to
\begin{equation}  \label{E:equivround}
   d_A^{\tp} - \tmu d_A^{\tm} : T_A \ra \Om^{\tp}.
\end{equation}

\begin{prop}   \label{P:invround}
The operator \eqref{E:wrapround} is invertible for $p$ close enough to $2$.
\end{prop}
\begin{proof}
It is enough to show that \eqref{E:equivround} is invertible.  Notice that
$(1-\tmu S) : \tL^p \ra \tL^p$ is invertible for $p$ close enough to $2$, by
lemma \ref{L:isomLp} and the estimate \eqref{E:normmu} on $\mu$. But an inverse
for \eqref{E:equivround} is given by $Q (1-\tmu S)^{-1}$ (see \cite{KM})
and we are done.
\end{proof}

Consider now the wrapped-up operator with the weighted adjoint of $d_A$
\begin{equation}  \label{E:wwrap}
 D_A^\de = d_A^{\tst_{\de}} + d_A^+ : \tL^p_{1,\de} (\Om^1(\WS ; \fg_E)) \ra 
      \tL^p_{0,\de} (\Om^0(\WS ; \fg_E)) \oplus \tL^p_{0,\de} (\Om^+(\WS ; \fg_E))
\end{equation}

\begin{prop}  \label{P:invirredw}
Let $A$ be an acyclic flat connection. Then, for $| \de |$ small enough and $p$ close enough
to $2$, the operator
$D_A^\de$ defined by \eqref{E:wwrap} is invertible.
\end{prop}
\begin{proof}
Since the set of invertible operators is open, this follows from proposition
\ref{P:invround}.
\end{proof}

We now assume that the connection $A$ is abelian. In that case the problem
splits globally:
\begin{equation}  \label{E:splitabel}
  D_A^\de = d_A^{\tst_{\de}} + d_A^+  \, = \,
   \begin{cases} 
      d^{*_{\de}} + d^+  \; \qquad \text{on the diagonal component of $\fg_E$}, \\
      d_{2\a}^{\cst_{\de}} + d_{2\a}^+ \qquad 
          \text{on the off-diagonal component of $\fg_E$},
  \end{cases}
\end{equation}
where $d_{2\a}$ is the exterior covariant derivative associated to the 
abelian connection with meridinal holonomy equal to $e^{-4i\pi \a}$.

\begin{lemma}  \label{L:diagabel}
The operator
$$
    d^{*_{\de}} + d^+ : L^p_{1,\de} (W) \ra L^p_{\de} (W)
$$
is invertible for $|  \de | \neq 0$ and small enough.
\end{lemma}
\begin{proof}
This follows from the standard techniques of \cite{LM}, see \cite{Floer}.
\end{proof}

\begin{lemma}  \label{L:invoffabel}
The operator
$$
  d_{2\a}^{\cst_{\de}} + d_{2\a}^+ : W^p_{1,\de} (\WS) \ra L^p_{\de} (\WS)
$$
is invertible for $| \de |$ small enough and $p$ close enough to $2$.
\end{lemma}
\begin{proof}
Proceeding as in the proof of proposition \ref{P:invLapcyl}, we can show
that $\check{\D}^W_{2\a,k}$ (the off-diagonal component of $\tDW_{2\a,k}$) is
invertible.  Indeed, $\check{\D}^M_{2\a,k}$ is invertible, and thus positive, 
because identifications similar to the ones obtained in  the
proof of proposition \ref{P:hommixed} hold.  Therefore, 
$$
  d_{2\a}^{\cst} + d_{2\a}^{\check{+}} : W^p_1(\WS) \ra L^p(\WS)
$$
is invertible. Mimicking the proofs of lemmata
\ref{L:isometry} , \ref{L:isomLp} and proposition \ref{P:invround}, shows that
$$
  d_{2\a}^{\cst} + d_{2\a}^{+} : W^p_1(\WS) \ra L^p(\WS)
$$
is invertible.  The lemma follows by the perturbation argument of proposition
\ref{P:invirredw}.
\end{proof}

Combining lemmata \ref{L:diagabel} and \ref{L:invoffabel}, we get 
\begin{prop}  \label{P:invabel}
Let  $A$ be the abelian connection.  For $| \de | \neq 0$ small enough and $p$ close
enough to $2$, the operator
\begin{equation}  \label{E:opw}
  D_A^\de = d_A^{\tst_{\de}} + d_A^+ : \tL^p_{1,\de} (\WS) \ra \tL^p_{0,\de} (\WS)
\end{equation}
is invertible.
\end{prop}

%
%
%
%


\section{Global theory.}  \label{SS:globalcyl}

We are now going back to the general set-up described in section \ref{SS:setupcyl}.

\subsection{The slice theorem.}
Let
\begin{equation}  \label{E:wadjcomp}
   d_{A^\a}^{\tst_\de} = w_{\de}^2 d_{A^\a}^{\tst} w_{\de}^{-2} : \Om^1(\WS; \fg_E) \ra
 \Om^0(\WS; \fg_E).
\end{equation}
be the weighted adjoint of $ d_{A^\a}$ in the mixed inner product and let
\begin{equation}   \label{E:wLapcomp}
  \D^W_{A^\a,0,\de} = d_{A^\a}^{\tst_{\de}}  d_{A^\a}
\end{equation}
be the corresponding Laplacian on sections of $\fg_E$.

\begin{prop}  \label{P:LapFredcyl}
Assume $\de_j \neq 0$ if $A_j$ is abelian. Then, the operator
$$
 \D^W_{A^\a,0,\de} : \tL^p_{2,\de} (\WS; \fg_E) \ra \tL^p_{0, \de} (\WS ; \fg_E)
$$
is Fredholm and self-adjoint.
\end{prop}
\begin{proof}
The proof is the same as before, gluing local parametrices on the compact part of $W$
and on the cylindrical ends.  Parametrices on
the compact part can be obtained as in chapter \ref{S:slice} because the normal
operator for the off-diagonal component of $r^2 \D^W_{A^\a,0,\de}$ is the same
as in the unweighted case.  Local inverses on the ends were obtained in
proposition \ref{P:invwLap}. 
Self-adjointness follows from elliptic regularity.
\end{proof}

Let us now consider a general singular connection $A = A^\a + a \in \mathcal{A}^{\a,p}_{\de}$. 
Let 
\begin{equation}  \label{E:fLapcyl}
  \tD_A : \tL^p_{2,\de} \ra \tL^p_{0,\de}
\end{equation}
be the Laplacian on sections associated to $A$ in the mixed metric.
\begin{prop}   \label{P:genLapcylFred}
Assume that $| \de |$ is small enough, with $\de_j \neq 0$ if $A_j$ is abelian.  Then, for
$p$ close enough to
$2$, the operator
\eqref{E:fLapcyl}is Fredholm and self-adjoint.
\end{prop}
\begin{proof}
Starting from proposition \ref{P:LapFredcyl}, the standard perturbation
argument of
\cite{TaubesCasson} based on the the openness of the space of Fredholm operators can be
applied.
\end{proof}

The proof of the following goes as before.
\begin{cor}  \label{C:Banachcyl}
With $\de$ as in proposition \ref{P:LapFredcyl}, the space
$(\tcB^{\a,p}_{\de})^*$ is a smooth  Hausdorff Banach manifold.  A slice for the
action of $\tcG^p_{\de}$ at $ A \in \mathcal{A}^{\a,p}_{\de}$ is given by
$\ker \, d_A^{\tst_\de}$.
\end{cor}

Combined with \eqref{E:goodconn}, we get
\begin{cor}
With $\de$ as in proposition \ref{P:LapFredcyl}, the space
$(\Bap_{\de})^*$ is a smooth  Hausdorff Banach manifold.
\end{cor}

\subsection{The deformation complex and the moduli spaces of ASD connections.}
As usual the deformation complex for $\tM^{\a,p}$ is 
\begin{equation}  \label{defcompcyl}
 0 \ra \tL^p_{2,\de}( \WS ; \fg_E) \overset{d_A}{\lra}
     \tL^p_{1,\de}( \Om^1( \WS ; \fg_E))   \overset{d_A^+}{\lra}
    \tL^p_{0,\de}( \Om^+( \WS ; \fg_E)) \ra 0.
\end{equation}
In view of corollary \ref{C:Banachcyl}, the deformation complex can be wrapped up
into the operator
\begin{equation}  \label{E:fwrapcyl}
  D_A^\de = d_A^{\tst_\de} + d_A^+ : \tL^p_{1,\de} \ra \tL^p_{0,\de}
\end{equation}

We first examine the operator \eqref{E:fwrapcyl} for the background connection $A^\a$.
\begin{prop}  \label{P:Fredcyl}
For $p$ close enough to $2$,  $|\de_j|$ small enough and $\de_j \neq 0$ if 
$A_j$ is abelian, the operator
\begin{equation}  \label{E:opwcomp}
  D_{A^\a}^\de = d_{A^\a}^{\tst_{\de}} + d_{A^\a}^+ :
     \tL^p_{1,\de} (\WS) \ra \tL^p_{0, \de}(\WS)
\end{equation}
is Fredholm.  
\end{prop}
\begin{proof}
In the usual way, one can glue local parametrices for \eqref{E:opwcomp} obtained on
the compact part of $W$ by proposition \ref{P:ASDwrapup} and on the cylindrical ends
by propositions \ref{P:invirredw} and \ref{P:invabel} to get a right parametrix for
$D_{A^\a}^\de$.  This shows that $D_{A^\a}^\de$ has a closed range of finite codimension.
A right parametrix for the adjoint of \eqref{E:opwcomp} can be obtained by duality
from the local parametrices. It follows that $\ker \, D_{A^\a}^\de$
is finite-dimensional. 
\end{proof}

The general case follows as in the proof of proposition
\ref{P:genLapcylFred}.
\begin{prop} For $p$ and $\de$ as in proposition
\ref{P:Fredcyl}, the operator \eqref{E:fwrapcyl} is Fredholm.  
\end{prop}

We now have all the ingredients to analyze the moduli spaces of
singular ASD connections on $(W, \S)$.
\begin{prop}  \label{P:modASDcyl}
Let $U$ be a precompact open set of $\WS$.  Then, for a generic metric perturbation
on $U$, the moduli space $(\tM^{\a,p}_\de)^*(W, \S)$ is a smooth finite dimensional
manifold, whose dimension is given by the index of $\Dad$.
\end{prop}
Proposition \ref{P:invindex} will show that the index of $\Dad$  doesn't depend
on $p$ or $\de$, when $p$ is close enough to $2$ and the $\de_j \geq 0$ are small enough,
with strict inequality if $A_j$ is abelian.
\begin{proof}
The only point worth commenting on is the smoothness of the moduli space. By duality,
$$
 \coker \, d_A^+ = \{ \om \in \tL^q_{0,-\de}( \Om^+) \, \mid \, d_A \om = d_A^*\om = 0 \},
$$
which shows that the elements in the cokernel of $d_A^+$ are harmonic forms $\om$. So, by 
Aronszajn's theorem, $\om$ vanish  on an open set, unless $\om \equiv 0$. Therefore,
the usual line of argument from \cite{FU} or \cite{DK} carries through.
\end{proof}

Combining proposition \ref{P:modASDcyl} with \eqref{E:gASD}, we get 
\begin{thm}   \label{T:ASDcyl}
Let $U$ be a precompact open set of $\WS$.  Then, for a generic metric perturbation
on $U$, the moduli space $(M^{\a,p}_\de)^* (W, \S)$ of irreducible singular ASD
connections is a smooth finite dimensional manifold, whose dimension is given by the index
of $\Dad$, decreased by the number of abelian limiting connections.
\end{thm}


%
%

\section{Some index formulae.}  \label{SS:index}

In this section, we first prove the invariance of the
index of $D^\de_{A^\a}$ in $p$ and $\de$, and thus that the moduli
spaces of ASD connections have a well-defined dimension.  After that, we
show how the index changes when two ends are glued together, or when one applies
a gauge transformation to a limiting connection on some end.  Those results
lead to defining for singular connections an analogue to Floer's grading
function.  This grading function exhibits a periodicity of $4$, instead 
of $8$ in the non--singular case \cite {Floer}.

\subsection{Invariance of the index.}
The first point that we would like to address is the
invariance of the index of
$D_{A^\a}^\de = \tdad + \dap$.  

\begin{prop}   \label{P:invindex}
For $p$ close enough to $2$ and $ | \de |$  small enough, with $\de_j \neq 0$ if
$A_j$ is abelian, the index of
\begin{equation}  \label{E:regop}
 \Dad = \tdad +\dap : \tL^p_{1,\de}( \Om^1(\WS ; \fg_E)) \ra
      \tL^p(\Om^0(\WS; \fg_E)) \oplus \tL^p(\Om^+(\WS; \fg_E))
\end{equation}
is independent of $p$.  It is also independent of $\de$ unless some limiting
connection $A_j$ is
abelian and the sign of $\de_j$ changes.
\end{prop}

We will determine later how the index changes when the weight corresponding to an
abelian connection changes sign.

\begin{proof} In fact,
$\tdad +\dap$ is homotopic  to
$\tda +\datp$ through Fredholm operators by progressively changing the metric from the
cone-angle
$\ak$ to the round metric on the off-diagonal component of $\fg_E$.  This is similar to
an argument of \cite{KM}.  Therefore, we simply need to study the index of $\tdad +
\datp$.  By continuity, we can assume that $\de_j = 0$ if $A_j$ is irreducible.  

Composing $\tda + \datp$ with its adjoint on the left, we have
that $\ker  (\tda + \datp) \subseteq \ker \, \tD_{A^\a,1}$.    
Let $p' \neq p$.  Elliptic regularity shows that if $u \in \ker \, \tD_{A,1} \cap
\tL^p_1( \WS)$, then $u$ is in $\tL^{p'}_1$ on the compact part of $W$ and on all
the strips $B_n = [n-1, n+1] \times M_j, n = 1, 2, \dots$  On each end $\MjR$,
$\tD_{A^\a,1}$ is given by
\eqref{E:Lapcylinv}, whence an eigenvalue decomposition shows that $u$ decays
exponentially. Therefore $u \in \tL^{p'}_1( \WS)$.

By duality the cokernel of  $\tda + \datp$ is the kernel of
\begin{equation}  \label{E:regcoker}
d_{A^\a} + \tda : \tL^q(\Om^0(\WS; \fg_E)) \oplus \tL^q(\Om^{\tp}(\WS; \fg_E))
     \ra \tL^q_{-1}( \Om^1(\WS ; \fg_E)).
\end{equation}
Composing \eqref{E:regcoker} with its adjoint on the left, the cokernel is
in the kernel of
\begin{multline}  \label{E:superLap}
  \begin{pmatrix}  \tD_{A^\a,0}  &(F_{A^\a}^{\tp})^{\tst}  \\
             F_{A^\a}^{\tp}  &\tD_{A^\a,2} \end{pmatrix} \, :
   \tL^q(\Om^0(\WS; \fg_E)) \oplus \tL^q(\Om^{\tp}(\WS; \fg_E))  \\
      \lra \tL^q_{-2}(\Om^0(\WS; \fg_E)) \oplus \tL^q_{-2}(\Om^{\tp}(\WS; \fg_E)).
\end{multline}
The off-diagonal operators of \eqref{E:superLap} are smooth zeroth order multiplication
operators,  by construction of
$A^\a$.  Therefore, all the elliptic regularity arguments applied before carry through for
the compact part of
$W$ and the strips $B_n$.  Moreover, these off-diagonal components also vanish on the
cylindrical ends since
$A^\a$ is flat there.  Hence, the exponential decay argument just
described can also be applied.  Therefore the cokernel of $\tda + \datp$
is also independent of $p$.

When $A_j$ is abelian, \eqref{E:regop} splits as in \eqref{E:splitabel}.  The
off-diagonal component can be homotoped until the weight corresponds to $\de=0$ and we
can  reason as above.  For the diagonal component, we can use for instance 
the standard theory of \cite{LM} to show that the kernel and cokernel don't depend on 
$p$ and can jump only when $\de=0$. It is also possible to adapt the argument given
above.
\end{proof}

In view of proposition \ref{P:invindex}, we can define the \emph{index} of an adapted
bundle pair to be the index of $\Dad$ when all the weights $\de_j$ are small positive
numbers.

\subsection{Index additivity.}
We are now studying the problem of gluing cylindrical ends. 
Let $(E,L) \ra (W, \S)$ be an adapted bundle with background connection $A^\a$ over a
manifold with cylindrical ends which is not necessarily connected. Assume that $W$ has an
end modelled on $(M,K) \times
\R^+$ and another end modelled on $(\bar{M}, \bar{K}) \times \R^+$, where $(\bar{M},
\bar{K})$ is diffeomorphic to
$(M,K)$, but has the opposite orientations.  Let $T>0$.  We can then form a new
manifold $(W_T^\sh, \S_T^\sh)$ by cutting off the ends at $(M,K) \times \{ T \}$ and
$(\bar{M}, \bar{K}) \times \{ T \}$ respectively and gluing $(M,K) \times \{ T \}$ and
$(\bar{M}, \bar{K}) \times \{ T \}$ together.  Assume that the flat connections
on the two ends agree and denote that common
limiting flat connection by $A_M$. We can then form a new adapted bundle
pair
$(E_T^\sh, L_T^\sh)
\ra (W_T^\sh, \S_T^\sh)$, with background connection $(A^\a_T)^\sh$. Let $(D_T^\de)^\sh$ be
the operator \eqref{E:opwcomp} associated to $(A^\a_T)^\sh$. The index of the bundle pairs
$(E_T^\sh, L_T^\sh )$  and $(E,L)$ are related as follows :

\begin{thm}  \label{T:indexadd}
If $A_M$ is irreducible, then ind$(E,L)$ = ind$(E_T^\sh, L_T^\sh )$. If $A_M$ is abelian,
then ind$(E_T^\sh, L_T^\sh )$ = ind $(E,L) - 1$.
\end{thm}

\begin{proof}
For the ease of exposition, identify $(\bar{M}, \bar{K}) \times \R^+$ with
$(M,K) \times \R^-$. Let $t$ be the axial coordinate.  Let $B_M$ be the operator
\begin{equation}  
 B_M = \begin{pmatrix} 0 &d^{\tst}_{A_M}   \\
          d_{A_M} &*d_{A_M} \end{pmatrix} :
   \Om^0( \MK) \oplus \Om^1( \MK) \ra \Om^0( \MK) \oplus \Om^1( \MK) 
\end{equation}

Let us first examine the case when
$A_M$ is irreducible. In that case, we can assume that $\de_M$, the weight corresponding to
the two ends  in question, is zero.  Then, making the usual identifications on both ends,
$\Dad$ looks like
\begin{equation}  \label{E:Dend}
  \Dad = \d_t + B_M.
\end{equation}
The proof is then perfectly standard \cite{DFK}.  We recall the main steps.  First,
after stabilization, we can assume that $D_{A^\a}$ is surjective.  For $T$ large enough, we
can then construct an approximate right inverse to $D_T^\sh$, showing that it is also
surjective. Finally, using well-chosen cut-off functions and the approximate right inverse,
it is possible to construct two approximately isometric injections $\ker \, D_T^\sh \ra \ker
\, D_{A^\a}$ and $\ker \, D_{A^\a} \ra \ker \, D_T^\sh$.

This proof doesn't work when $A_M$ is abelian, because we cannot take $\de_M$ to be zero.
Instead $\de_M$ must be a small positive number.  It follows that the identification
\eqref{E:Dend} doesn't hold over the ends.  Over one end,
\begin{equation}  \label{E:Dend1}
  \Dad = \begin{pmatrix}  \de_M &0 \\  0 &0 \end{pmatrix} + \d_t + B_M, 
\end{equation}
while over the other end,
\begin{equation}  \label{E:Dend2}
  \Dad = \begin{pmatrix} -\de_M &0 \\  0 &0 \end{pmatrix} + \d_t + B_M. 
\end{equation}
So we need to introduce a bridge $(M,K) \times \R$ between the two ends and to compute the
index of the operator 
\begin{equation}  \label{E:bridge}
  D_{A_M}^{-(\de_M, \de_M)} : \tL^p_{1,(\de_M, \de_M)} \ra  \tL^p_{0,(\de_M, \de_M)}
\end{equation}
which agrees with \eqref{E:Dend1} over $M \times (-\infty,-1)$ and with \eqref{E:Dend2}
over $M \times (1, \infty)$.  As usual the operator \eqref{E:bridge} splits globally
between on-- and off--diagonal components, as in \eqref{E:splitabel}.  The off-diagonal
component can be homotoped until $\de_M = 0$, showing that its index is zero.  The
theory of \cite{LM} shows that the index of the diagonal part
is $-1$.  Hence, the index of the bridge operator is $-1$.  We can then apply twice
the proof sketched in the irreducible case to obtain the index formula.
\end{proof}

\subsection{A change of gauge formula.}
From theorem \ref{T:indexadd}, we can also deduce a
change of gauge formula.
\begin{thm}  \label{T:gaugeindex}
Let $(W,\S)$ be a manifold with cylindrical ends, one of the ends being $(M,K) \times
\R^+$.  Let $A^\a$ be a background connection for an adapted bundle over $(W, \S)$, the
limiting connection on $(M,K)$ being denoted $A_M$. Let
$g$ be a gauge transformation in $\mathcal{G}(M,K)$. Lift $g$ to $(M,K) \times \R^+$ and extend
it to a gauge transformation $\tg$ on $(W,\S)$, arranging that $\tg$ be the identity on the
other ends. Then $\tg^* \! A^\a$ is a new background connection, and
\begin{equation}   \label{E:gaugeindex}
  \text{ind} \, D^\de_{\tg^* \! A^\a} = \text{ind} \,  D^\de_{A^\a} + 8 \, \deg \, g
       - 4 \, \deg \, g|_K.
\end{equation}
\end{thm}
\begin{proof}
In view of theorem \ref{T:indexadd}, it is enough to prove the special case where $(W, \S) =
(M,K) \times \R$ and $A^\a = A_M$.  To prove that case, we apply theorem \ref{T:indexadd}
again, gluing the two ends of the cylinder to obtain a pair of compact manifolds
diffeomorphic to $(M~\times~S^1, K ~\times~S^1)$.  The characteristic numbers of the bundle
pair obtained by gluing with the gauge transformation $g$ are given by
$$
  k = \deg \, g  , \qquad \qquad  \ell = - \deg \, g|_K.
$$
Therefore, \eqref{E:gaugeindex} follows from formula \eqref{E:dimASD} giving
the index of $D_{A^\a}$ in the compact case.
\end{proof}

\subsection{The grading function.}

Let us go back to the special of the cylinder $(M,K) \times \R$.  For $A,B$ two 
isolated flat connections on $(M,K)$, consider the adapted bundle pair $(E,L)$ over
$(M,K) \times \R$, with $A$ as the limiting connection at $-\infty$ and $B$ as the
limiting connection at $+\infty$.  Let $\tmu (A,B)$ be the formal dimension of
the correponding moduli space of ASD singular connections, that is
$$
 \tmu (A,B) = \text{ind} (E,L) - \dim \G_A - \dim\G_B ,
$$
where $\G_A$ and $\G_B$ are the isotropy subgroups of $A$ and $B$ respectively.
Theorem \ref{T:indexadd} can be rephrased as
\begin{equation}    \label{E:muadd}
 \tmu(A,C) = \tmu(A,B) + \tmu(B,C) + \dim \G_B ,
\end{equation}
which has the same form as the well-known addition formula for the spectral flow in
Floer homology. It follows from \eqref{E:muadd} that
\begin{equation}    \label{E:muAA}
 \tmu(A,A) = - \dim \G_A.
\end{equation}
If $\TH$ is a representative of the abelian connection, we can define a grading
function
\begin{equation}  \label{E:tmdef}
 \tmu(A) = \tmu( \TH, A) + \dim \G_A.
\end{equation}
It follows from \eqref{E:muAA} and \eqref{E:muadd} that 
$$
 \tmu( \TH ) = 0
$$
and
$$
 \tmu(A,B) = \tmu(B) - \tmu(A) - \dim \G_B.
$$
To obtain the $\Z / 4$ grading function, we just need to observe that
theorem \ref{T:gaugeindex} can be rephrased as
\begin{equation}  \label{E:mugauge}
 \tmu (g^{*} A) = \tmu (A) + 8 \, \deg \, g - 4 \, \deg \, g|_K
\end{equation}
if $g$ is a gauge transfomation on $(M,K)$. In view of \eqref{E:mugauge}, when
dividing by the action of the gauge group, we obtain a $\Z / 4$ grading
function $\mu : \mathcal{F}^\a \ra \Z / 4$ by
$$
 \mu([A]) = \tmu (A) \text{mod} \, 4.
$$
The previous discussion shows that this function displays the properties
of a spectral flow.
\vfil\eject
\chapter*{Loose ends and speculations.}

In this afterword, we discuss a number of points which should be addressed in
order to complete the program of developing a Floer homology for singular
connections.  The main missing element is an analysis of the ends of the moduli
spaces of singular ASD connections on the cylinder.  Based on theorem \ref{T:indexadd},
it is possible to prove a gluing theorem for suitable families of singular ASD
connections, just as in the usual case \cite{DFK}.  It is also possible to prove
that one obtains in that way local diffeomorphisms between the corresponding
moduli spaces.  However, there are a couple of technical obstacles to proving that
the ends of the moduli spaces are as one would expect them to be.  One of those
obstacles is that we don't have uniform decay estimates for singular ASD
connections on the cylinder.  The problem comes indeed from the fact that we
couldn't prove enough regularity for singular ASD connections.  The other
obstacle comes from the perturbations.  In order to get the usual cusp--trajectory
picture for Uhlenbeck compactness in Floer homology, one needs to have smooth
moduli spaces for a translation--invariant equation on the cylinder.  Perturbing
the metric as in theorem \ref{T:ASDcyl} cannot always be done in a 
translation--invariant way.  Instead, one needs to consider holonomy perturbations.
The proof of smoothness of the moduli spaces of perturbed ASD connections relies
on a gauge--theoretic version of Aronszajn's unique continuation theorem for
the Cauchy problem, as proved in \cite{TaubesAron}.  We couldn't adapt the proof
given there, mainly again because of the lack of regularity of the singular ASD
connections.  Another problem arising from perturbations is that we couldn't show
the invertibility of the operator $D_A$ defined in \eqref{E:wrapround} in
the translation--invariant case (the model case on the ends) if $A$ is a
perturbed flat connection.  More precisely, the problem comes from lemma
\ref{L:isometry}, which doesn't hold if $A$ is not flat.  Therefore, the 
conformal structure for the off--diagonal component cannot be changed back
to the round metric.  However, it is possible that the operator $D_A$ could
be analyzed by a generalization of edge theory to manifolds with corners,
using a hyperbolic metric on the off--diagonal component as in \cite{R1}.

We should mention that in his thesis \cite{C}, Olivier Collin developed an
instanton Floer
homology for three--manifolds with an orbifold singularity along a knot.
Our theory would likely be equivalent to his for a rational value of the
holonomy parameter $\a$.  The advantage of ours is that it would allow
to interpolate between different rational values of $\a$, possibly proving
with a cobordism argument that the Floer homology is invariant in $\a$ as 
long as no square root of the Alexander polynomial is crossed.  It seems
likely that the Euler characteristc of the Floer homology in question would
be Herald's Casson invariant for knots \cite{Herald2}.   Li \cite{Li} has a symplectic
Floer homology for knots correponding to rational values of $\a$.  We can formulate
a conjecture \`{a} la Atiyah--Floer stating that Li's theory should be isomorphic
to Collin's and/or ours.

We would like to conclude by suggesting a possible application of the theory to
obtain bounds on the unknotting number of a knot.  The unknotting number is the
minimal number of crossing changes to apply to a planar projection of a knot in
order to obtain the unknot.  Crossing changes can be performes by a $\pm 1$ surgery
on a small trivial knot linking the crossing.  Floer \cite{Floer2} studied the
effect of surgery on his homology and encoded it in an exact triangle (see \cite{BD}
for details).  If his techniques can be adapted to singular connections, the
singular Floer homology might provide information on the unknotting number.

\backmatter

\vfil\eject

\appendix
\chapter{Elliptic edge operators in a nutshell.} \label{A:edge}

\setcounter{equation}{0}
\renewcommand{\theequation}{A.\arabic{equation}}

This appendix contains a summary of the theory of
 elliptic edge operators, as developed in \cite{Mazzeo}.  It also corrects
a few misprints from that paper, and contains a proof of the key fact
 (proposition \ref{P:Lp}) necessary to extend
the theory to $L^p$ spaces for $p \ne 2$. The last section develops elliptic
estimates for operators which become elliptic edge only after some
alteration.  It is this kind of operators which appear in the body of the thesis.

Let $M$ be a Riemannian compact manifold with boundary $\d M$.  Assume that $\d M$ is the
total space of a locally trivial fibration
\begin{align}
    F  \longrightarrow & \, \d M   \notag \\
                       & \, \, \downarrow \notag \\
                       & \, \,  B  \notag
\end{align}
with fibre $F$ and base space $B$.  The algebra of edge differential operators
is defined as the algebra of operators generated over $C^{\infty}(M)$ by the space
$\Vedge$ of
smooth vector fields on M which are tangent to the fibres of $\d M$.  More
concretely, let $x$ be a defining function for $\d M$, $y$ be local coordinates
on $B$ lifted to $\d M$ and $z$ be local coordinates along the fibers of $\d M$.
Extend $y$ and $z$ locally to the interior of M.  Then a differential operator
 $L$ of order $m$ is \emph{edge} if, locally near $\d M$, it takes the form
$$
  L \quad = \sum_{j+| \a | + | \b | \, \leq \, m} a_{j, \a , \b  } (x,y,z)
            (x \d _x)^j (x \d _y)^{\a} \d _{z}^{\, \b} ,
$$
where the $a_{j, \a , \b  }$'s are smooth functions, possibly matrix-valued when
dealing with operators on sections of vector bundles.  The operator $L$ is said
to be {\em elliptic \/}  if its \emph{edge principal symbol}
$$
{}^e\sigma_m
(L)(x,y,z)(\x , \h , \z ) \quad = \sum_{j+| \a | + | \b | \, = \, m} a_{j, \a ,
\b  } (x,y,z)
              (i \x)^j (i \h)^{\a} (i \z)^{\b}
$$
is elliptic, i.e. ${}^e\sigma_m
(L)(x,y,z)(\x , \h , \z)$ is invertible for
$(\x , \h , \z ) \neq 0$.

To an edge operator, we can associate a family of indicial operators, 
parametrized by $y \in B$, and defined over $F_y \times \R ^{+}$.  The family is given
in local coordinates by
\begin{equation}   \label{E:indicial}
   I(L) \quad = \sum_{j+| \b | \, \leq \, m} a_{j,0,\b}(0,y,z)
            (s \d_s)^j \d_z^{\, \b},
\end{equation}
where $s \in (0, +\infty)$ is the $\R^{+}$ coordinate.  We can associate to
the family of indicial operators, its indicial family (basically its
conjugate by the Mellin transform in the $s$-direction), parametrized by $(\z, y) \in \C \times B$, and defined
over $F_y$
\begin{equation}  \label{E:indMellin}
   I_{\z}(L) \quad =\sum_{j+| \b | \, \leq \, m} a_{j,0,\b}(0,y,z)
            \z^j \d_z^{\, \b}.
\end{equation}
Notice that $I_{\z}(L)$ is elliptic on $F_y$ if $L$ is elliptic.  The 
{\em boundary spectrum \/} of $L$ is defined as
$$
  spec_b(L) = \{ \, \z \in \C \, \mid  I_{\z}(L) \text{ is not invertible on } L^2(F_y)
          \text{ for some } y \in B \}.
$$
In the cases studied in this thesis, all the operators in the indicial family are
similar;  therefore, $spec_b(L)$ forms a discrete set in $\C \,$,
with $\Re \z \rightarrow \infty$ if $\Im \z \rightarrow \infty$. 
  We also define $\Re spec_b(L)$ to be the discrete subset of $\R$ obtained by
the projection of
$\sb (L)$ on the real axis.  From the theory in  \cite{Melrose}, it follows that
$I(L)$ is invertible on $s^{\de}L^2(F_y \times \R ^{+})$ iff $\de \notin \rsb (L)+1/2$.

There  is another family of operators associated to an edge operator,
 the family of {\em normal operators\/} $N(L)$.  It is parametrized by $y \in B$ 
and is defined over $F_y \, \times \R^+ \times \R^k$, where $k= \dim B$.
\begin{equation}  \label{E:normal}
  N(L)  \quad = \sum_{j+| \a | + | \b | \, \leq \, m} a_{j, \a , \b  } (0,y,z)
            (s \d_s)^j (s\d_u)^{\a} \d_{z}^{\, \b} ,
\end{equation}
where $s$ is the coordinate on $\R^+$, $u$ denotes the coordinates on $\R^k$,
and $y$ is just the parameter in $B$. Again, in our cases, all the normal operators
are similar.  Using the homogeneity properties of $N(L)$
in $(s,u)$, the invertibility properties of $N(L)$ follow from the
invertibility properties of what is essentially its conjugation by the Fourier transform
in the $u$ variables,
\begin{equation}   \label{E:normalFourier}
  \widehat{N(L)}  \quad = \sum_{j+| \a | + | \b | \, \leq \, m} a_{j, \a , \b  } (z)
            (s \d_s)^j (is \hat{\h})^{\a} \d_{z}^{\, \b},
\end{equation}
where $\hat{\h}$ is a parameter in $S^{k-1}$, the unit sphere in $\R^k$.  In our cases,
$N(L)$ will be isotropic in $u \in \R^k$, so that it is enough to consider
the case
$\hat{\h}=(1,0,\dots,0)$, which yields a \emph{Bessel-type} operator
\begin{equation}   \label{E:Besseltype}
    L_0 \quad = \sum_{j+| \a | + | \b | \, \leq \, m} a_{j, \a , \b  } (z)
            (t \d_t)^j (it)^{\a_1} \d_{z}^{\, \b}.
\end{equation}
The operator $L_0$ is  Fredholm
over  $t^{\de}L^2(F \times \R ^{+})$ if $\de \notin \rsb (L)+1/2$.  Moreover,
its kernel and cokernel don't change unless $\de$ crosses an element of
 $\rsb (L)~+~1/2$.  Therefore, there is a (possibly empty) maximal interval 
$( \underline{\de} \, , \overline{\de} )$ so that $L_0$ is invertible on
$t^{\de}L^2(F~\times~\R^{+})$ for 
$\de~\in~( \underline{\de} \, , \overline{\de} ), \, \de~\notin~\rsb (L) + \half.$  Note that
$\underline{\de} \, , \overline{\de} \in \rsb (L) + \half$.

Before we talk about spaces of pseudodifferential operators on $M$, we need to introduce 
the notion of asymptotic development.   Let $E$ be a discrete subset of $\C$ such that
$\Re z \rightarrow \infty$ if $\Im z \rightarrow \infty \, , z \in E$. We call $E$ an 
{\em index set\/} and set 
$$
\inf E = \inf \{ \, \Re z \, , \: z \in E \}.
$$
Let $X$ be a 
manifold with boundary $\d X$ and $\r$ a defining function for $\d X$.  Then a function
$u \in \mathcal{A}^E_{phg}(X)$ if it admits a polyhomogeneous expansion given by $E$ at $\d X$:
\begin{equation}    \label{E:defpoly}
   u \quad \sim    \sum_{ \begin{smallmatrix}
                                z \, \in \, E, \\  \Re z \, \ra  \, \infty
                               \end{smallmatrix} }
                         \sum_{p=0}^{p=p_z}   \r^z (\log \r)^p \, a_{z,p} \, , \qquad
                a_{z,p} \in C^{\infty}(X)
\end{equation}
and $ \Re z \rightarrow \infty$ if $p_z \ra \infty$. This means that for all
$N > 0$, there is some $M$ so that 
$$ u \quad -  
         \sum_{ \begin{smallmatrix}
                    z \, \in \, E ,  \\
                    \Re z \, < \, M
                 \end{smallmatrix} } 
            \sum_{p=0}^{p=p_z}   \r^z (\log \r)^p \, a_{z,p}
$$
vanishes to order $N$ at $\d X$. The definition extends
 easily to sections of bundles.   For manifolds with corners, one chooses an index set
for each   face of codimension one, and iterates \eqref{E:defpoly}, 
see \cite{Melrose} .

Operators on $M$ have their kernels defined on the manifold with corners $M^2 = 
M \times M$.   However, the analysis is facilitated when those kernels are lifted
to another manifold with corners $M^2_e$, obtained by blowing up a certain submanifold
of the corner of $M^2$ (again, see \cite{Mazzeo} for  details).   $M^2_e$ has three
codimension one faces, denoted $M_{10}$, the {\em left face\/},
$M_{01}$, the {\em right face\/}, and $M_{11}$, the {\em front face\/}. The diagonal
 $\Delta$
of $M^2$ lifts as $\Delta_e$ to $M^2_e$, where it intersects $M_{11}$ transversely.
  There are three basic
spaces of edge pseudodifferential operators on M :
\begin{enumerate}
  \item  $\Psi^m_e(M)$, where $m \in \Z$. These are the operators of order $m$ in the 
         {\em small calculus\/}. Their kernels,
          considered on $M^2_e$,
          are smooth away from $\Delta_e$ and vanish to infinite order at $M_{10}$ and
          $M_{01}$.  They have a certain conormal singularity along $\Delta_e$ (see
          \cite{Mazzeo} for the exact description).  The kernel of an edge differential
          operator $L$ of order $m$ is supported on $\Delta_e$ and therefore, $L \in
         \Psi^m_e$.
  \item  $\Psi^{- \infty , \mathcal{E}}_e(M)$, where $\mathcal{E}=(E_{10},E_{11},E_{01})$ is a
         triplet of index sets.  The kernels
         of the operators in this space are  given by  functions which are smooth on the
         interior of $M^2_e$, with asymptotic developments at the face $M_{ij}$ corresponding
         to the index set $E_{ij}$.  We can then define
$$
    \Psi^{m , \mathcal{E}}_e(M) \quad = \quad \Psi^{- \infty , \mathcal{E}}_e(M) \,
       + \, \Psi^m_e(M).
$$
   \item $\Psi^{- \infty, \tilde{\mathcal{E}}}(M)$, where $\tilde{\mathcal{E}} = (E_{10},E_{01})$.
          The operators in this space are the {\em very residual \/} terms of
          the calculus and have smooth kernels on the interior of $M^2$,
          with asymptotic expansions at the left (right) face of $M^2$ given by 
          $E_{10}  \: (E_{01})$.  This allows us to define 
$$
     \Psi^{m , \mathcal{E}}(M) \quad = \quad \Psi^{m , \mathcal{E}}_e(M) \, 
         + \, \Psi^{- \infty, \tilde{\mathcal{E}}}(M).
$$
\end{enumerate}

For the composition properties of those spaces of operators, see \cite{Mazzeo}.

We now define the Sobolev spaces on which those operators act:
$$
  x^{\de}W^{p,k}_e(M) \quad = \{ u \in L^p(M) \mid u=x^{\de}v,\,  V_1 \dots V_k \,
       v \in L^p(M), \, V_i \in \Vedge (M)  \}.
$$
If $k$ is not a positive integer, the definition is obtained by duality and interpolation.
The extension to sections of vector bundles is straightforward, using local trivializations.
The dual of $x^{\de}W^{p,k}_e(M)$ is identified with $x^{-\de}W^{q,-k}_e(M), \; \frac{1}{p}+
\frac{1}{q} = 1$, under the $L^2$ inner product.
We also define $x^{\de}H^k_e = x^{\de}W^{2,k}_e$.  Now, a rescaling argument, as presented in
\cite{Mazzeo}, shows that
\begin{enumerate}
  \item  Operators in $\Psi^m_e$ map $x^{\de}W^{p,k}_e$ to $x^{\de}W^{p,k-m}_e$. In particular, this
         is true for edge differential operators of order $m$.
  \item  The corresponding statements of elliptic regularity hold for an elliptic edge 
         operator $L$ of order $m$, i.e. if $u,\,  Lu \in x^{\de}W^{p,k}_e$,
         then $u \in x^{\de}W^{p,k+m}_e$, with
         the  estimate:
\begin{equation}  \label{E:ellestim}
   \|u\|_{x^{\de}W^{p,k+m}_e}   \leq  C \, \|Lu\|_{x^{\de}W^{p,k}_e} + \|u\|_{x^{\de}W^{p,k}_e}.
\end{equation}
\end{enumerate}

Now, we denote by $L^t$ the formal adjoint of $L$ under the $L^2$ inner product, and by
$L^{*}=x^{2\de}L^t x^{-2\de}$ the adjoint in $x^{\de}L^2$. The boundary spectra of $L^t$
and $L^{*}$ can  be computed to be
\begin{equation}  \label{E:specadj}
 \begin{aligned}
     \sb (L^t) = & \; \{ -\bar{\z} - 1 \mid \, \z \in \sb (L) \}       \\
     \sb (L^{*}) = & \; \{ -\bar{\z} + 2\de - 1 \mid \, \z \in \sb (L) \}.   
 \end{aligned}
\end{equation}
This is \cite[formula (4.16)]{Mazzeo} where a sign mistake has been corrected
and the conjugations have been restored.

We are now ready to state the main result regarding elliptic edge operators
( \cite[ theorem (4.20)]{Mazzeo} (note that we corrected a mistake regarding
the index sets $F_{10}$ and $F_{01}$):

\begin{thm}  \label{T:Fred}
  Let $L$ be an elliptic differential edge operator of order $m$ with constant
indicial roots and isotropic normal operator.  Let
$\de~\in~( \underline{\de} \, , \overline{\de} )$,  $\de~\notin~\rsb (L)~+~1/2$, with
$\underline{\de} \, , \overline{\de}$ defined from $L_0$ as described 
above. Then the map 
$$
   L :  x^{\de}H^{\ell}_e \rightarrow x^{\de}H^{\ell -m}_e
$$
is Fredholm.  Denote by $P_1$ and $P_2$ the orthogonal projectors in $x^{\de}L^2$
onto the kernel and cokernel of L.  There is a generalized inverse $G$
for $L$, satisfying
\begin{equation}  \label{E:geninv}
      GL = I - P_1, \qquad  LG = I - P_2
\end{equation}
as operators on  $x^{\de}H^{\ell}_e$, for any $\ell \in \R$. We have   
\begin{gather}
 P_1 \in \Psi^{- \infty, \tilde{\mathcal{E}}}, \quad
 P_2 \in \Psi^{- \infty, \tilde{\mathcal{F}}},   \notag \\ 
 G \in \Psi^{- m, \mathcal{H}},\notag
\end{gather}
where
\begin{equation}     \label{E:asymptkernel}
 \begin{aligned}
	\inf E_{10} =  & \inf \, ( \rsb(L) \, \cap  \, ( \de - \tfrac{1}{2}, \infty) ),  \\   
 \inf E_{01} =  & \inf E_{10} - 2\de \, ,   \\
 \inf F_{10} =  & \inf \, ( \rsb(L^{*}) \, \cap \, ( \de - \tfrac{1}{2}, \infty) ),  \\   
 \inf F_{01} =  & \inf F_{10}  - 2\de \,  ,  \\
 \inf H_{10} =  & \min \{ \inf E_{10} , \inf F_{10} \},  \\
 \inf H_{01} =  & \min \{ \inf E_{01} , \inf F_{01} \},   \\
 \inf H_{11} =  & \, 0. 
 \end{aligned}
\end{equation}
The maps
\begin{align}
  G  \; &  : \;  x^{\de}H^{\ell}_e \rightarrow x^{\de}H^{\ell + m}_e  \notag \\
  P_i \; & : \; x^{\de}H^{\ell}_e \rightarrow x^{\de}H^s_e   \notag
\end{align}
are bounded for every $\ell, \, s \in \R$. Moreover, if $u \in x^{\de}L^2$ and 
$L u = 0$, then $u$ is polyhomogeneous with
$u \in \mathcal{A}^{E_{10}}_{phg}$.
\end{thm}

In order to obtain a good $L^p$ theory, $p \neq 2$, we need the following fact :
\begin{prop}   \label{P:GLp}
With notations as in theorem \ref{T:Fred}, let $\de'$ be such that
\begin{equation}    \label{E:ineqLp}
   \inf H_{10} > \de' - \frac{1}{p} \, ,  \quad  \inf H_{01} + \de' > - \frac{1}{q},
\end{equation}
with $\dfrac{1}{p}+\dfrac{1}{q} = 1$. Then the generalized inverse
$$
  G \: : \: x^{\de'}W^{p,k}_e \ra x^{\de'}W^{p,k+m}_e
$$
is well-defined and bounded and the equations \eqref{E:geninv} still hold. 
\end{prop}

Remark that in the statement of proposition \ref{P:GLp}, the operators
$G, P_1$ and $P_2$ and the index sets $H_{10}, H_{01}$ 
are the ones constructed with the weight $x^\de$, not
with the weight $x^{\de'}$.  We need to be able to change exponents 
 that way in the thesis.   
Proposition \ref{P:GLp} follows easily from the mapping properties of operators in the
small calculus $\Psi^m_e$ between
$L^p$ spaces obtained by rescaling, as explained above,
 and the following extension of \cite[Theorem (3.25)]{Mazzeo},
 as suggested there.
\begin{prop}   \label{P:Lp}
Let $A \in \Psi^{-\infty, \mathcal{E}}$ for some collection of index sets $\mathcal{E}$. Then
$$
      A  : x^{\de}W^{p,\ell}_e \ra x^{\de '}W^{p,s}_e
$$
is bounded for any $s, \ell \in \R$ provided
\begin{equation}  \label{E:kerasympt}
   \inf E_{10} > \de ' - \frac{1}{p},  \quad  \inf E_{01} + \de > - \frac{1}{q},
   \quad \de + \inf E_{11} \geq \de '.
\end{equation}
\end{prop}

\begin{proof}
As in \cite{Mazzeo}, to which we refer for notations, it is enough to consider the case of
operators acting on half-densities and we can reduce to the case
 where $s = \ell = 0$ and $\de = \de ' = 0$.
Hence, all we need to show is the estimate 
\begin{equation}  \label{E:estLp}
  | \langle A u,v \rangle |  \leq   C \|u\|_p \|v\|_q.
\end{equation}
This is straightforward if $A$ is very residual.  Therefore, we can assume that 
$A \in \Psi^{-\infty, \mathcal{E}}_e$ and consider its kernel $\k_A$ lifted to $M^2_e$.
The left hand side of \eqref{E:estLp} becomes
\begin{equation}  \label{E:estker} 
 | \int \k_A \, u_R \,  v_L \, \nu^2 | ,
\end{equation} 
where $\nu^2$ is a non-vanishing density on $M^2_e$, and the subscripts $R$ and $L$ indicate that
$u$ and $v$ have been lifted from the right and from the left, respectively (note there
 is a typo in the proof of \cite{Mazzeo}, the subscripts $R$
and $L$ having been interchanged).  If we denote by $\r_{10}$ and $\r_{01}$ defining
functions for the left and right faces of $M^2_e$, then \eqref{E:estker} can be estimated
as follows: 
\begin{equation}
\begin{split}
 \int | \k_A | \, | u_R | \, | v_L | \, \nu^2
   &= \int | \k_A | \; | u_R |   (\r_{01} / \r_{10})^{\frac{1}{p+q}} \; 
           | v_L  | ( \r_{10} / \r_{01} ) ^{\frac{1}{p+q}} \; \nu^2  \\
   &\leq  \left\{ \int | \k_A | \, | u_R |^p   ( \r_{01} / \r_{10} ) ^{\frac{1}{q}} \, \nu^2 \right\} ^{\frac{1}{p}}
       \left\{ \int | \k_A | \, | v_L  |^q (\r_{10} / \r_{01})^{\frac{1}{p}} \, \nu^2 \right\} ^{\frac{1}{q}} , \\
\intertext{using H\"{o}lder's inequality,} 
   &= \left\{ \int | \k_A | \, | u |^p   \r_{10}^{-\frac{1}{q}} \r_{01}^{-\frac{1}{p}} \, \g^2_R \, d\th \right\} ^{\frac{1}{p}}
      \left\{ \int | \k_A | \, | v |^q   \r_{10}^{-\frac{1}{q}} \r_{01}^{-\frac{1}{p}} \, \g^2_L \, d\th \right\} ^{\frac{1}{q}} , \\
\intertext{using \cite[equation (3.21)]{Mazzeo}  ($\g$ is a non-vanishing half-density on $M$
 and
$d \th$ is a non-vanishing density on $M_{11}$), } 
   &\leq   C \; \|u\|_p  \; \|v\|_q, \\
\end{split}  \notag
\end{equation}
by Fubini's theorem, given the asymptotics \eqref{E:kerasympt} of the kernel $\k_A$.
\end{proof}

It often happens, as in this thesis, that the operators studied become edge only after
multplying by some power of $x$.  For simplicity, let's assume we have a second order
 operator
$\D$ over the manifold $M$, such that $L = x^2 \D$ is elliptic edge.  In that case, the 
usual statement of elliptic regularity for elliptic edge operators \eqref{E:ellestim}
is not completely satisfactory.  Indeed, 
\begin{equation}   \label{E:mapD}
    \D \, : \, x^{\de+2} H^{k+2}_e  \ra x^\de H^k_e
\end{equation}
is bounded.  Therefore, what we would like is a statement of elliptic regularity that would
for instance read as follows : if $u, v \in x^\de H^k_e$ and $\D u = v$,
then $u \in x^{\de+2} H^{k+2}_e$,
with the corresponding estimates.
However, what we have in terms of the operator $L$ is
$$
  u \in x^\de H^k_e \quad \text{and} \quad 
     L u = x^2 v \, \in \, x^{\de+2} H^k_e \, \subset \, x^\de H^k_e \, , 
$$
and therefore, it can only be concluded from 
\eqref{E:ellestim} that $u \in x^\de H^{k+2}_e$.
Still, as we shall see,
it is possible to raise the power of $x$ from $\de$ to $\de+2$ as desired
if there is some elbow room in the distribution of the indicial roots of $L$.

\begin{prop}  \label{P:wellreg}
Let $\de$ satisfy the hypotheses of theorem \ref{T:Fred} for $L=x^2 \D$.  With the notations
of theorem \ref{T:Fred}, assume 
$\de'$ is such that
\begin{equation}    \label{E:wineqLp}
   \inf H_{10} > \de' + 2 - \frac{1}{p} \, ,  \quad  \inf H_{01} + 2 + \de' > - \frac{1}{q},
\end{equation}
If $u \in x^{\de'} W^{p,\ell}_e, \,  v = \D u \in x^{\de'} W^{p,k}_e$, then
$u \in x^{\de'+2} W^{p,k+2}_e$ and
\begin{equation}   \label{E:wellest}
   \|u\|_{x^{\de'+2} W^{p,k+2}_e} \,  \leq \, C \, \bigl( \,  
        \| \D u \|_{x^{\de'} W^{p,k}_e} +  \|u\|_{x^{\de'} W^{p,\ell}_e} \bigr) .
\end{equation}
\end{prop}

\begin{proof}
Let $w = x^2 v \in x^{\de'+2} W^{p,k}_e$.
 We apply $G$, the generalized inverse for $L$ corresponding to the
weight $\de$, to both sides of the equality    $L u = w$.  By proposition \ref{P:GLp},
$G w \in x^{\de'+2} W^{p,k+2}_e$. Moreover, we have from \eqref{E:geninv},
\begin{equation}  \label{E:GPu}
    u = Gw - P_1 u
\end{equation}
By proposition \ref{P:Lp}, $P_1 u \in x^{\de'+2} W^{p,k+2}_e$.  So the proposition, including
\eqref{E:wellest}, follows from \eqref{E:GPu}.
\end{proof}

From \eqref{E:asymptkernel} and \eqref{E:specadj}, the inequalities \eqref{E:wineqLp}
can be seen to be equivalent to saying that there is no indicial root
of $L$ whose real parts falls between $\de - \half$ and $\max \{ 2\de -\de' -2- \frac{1}{q},
\de' +2- \frac{1}{p} \}$, nor between
$\de - \half$ and $\min \{ 2\de - \de' -3 + \frac{1}{p}, \de' + 1 + \frac{1}{q} \}$.
This is what we meant by elbow room in the distribution of the indicial roots of $L$.

The only thing which is really used in the proof of proposition \ref{P:wellreg} is
the equality \eqref{E:GPu}, as well as the mapping properties of $G$ and $P_1$, which,
according to propositions \ref{P:GLp} and \ref{P:Lp}, depend exclusively on the distribution
of the indicial roots of $L$ around $\de - \half$.  It is therefore possible to obtain
other statements of elliptic regularity by the same kind of reasoning.  For instance,
we will need in the thesis:

\begin{prop}  \label{P:intellestim}
Let $\Om \subset \subset V$ be open subsets of $M$. If $u \in x^\de H^\ell_e(V) \, , \,
\D u \in x^{\de'}W^{p,k}_e (V)$, then $u \in x^{\de'+2}W^{p,k+2}_e (\Om)$, with
the estimate
\begin{equation}  \label{E:intellestim}
  \|u\|_{x^{\de'+2}W^{p,k+2}_e (\Om)} \, \leq \, C \, \bigl( \,
    \| \D u \|_{x^{\de'}W^{p,k}_e (V)} + \|u\|_{x^\de H^\ell_e(V)} \bigr) \, ,
\end{equation}
provided the elements of $\rsb(r^2 \D)$ are located far enough from $\de - \half$.
\end{prop}

We don't give the precise condition on the indicial roots here, since it is not needed in
the thesis.  In fact, what is shown there is that it is always possible to create as much
elbow room as desired by altering $\D$.  This proposition is really new
only when $\Om \cap \d M \neq \O$; otherwise, it's just a well-known form
of interior elliptic regularity.

\begin{proof}
The proof proceeds as in proposition \ref{P:wellreg}, using cut-off functions equal
to $1$ on $\Om$ and supported on $V$, as in one of the usual proofs of interior elliptic 
regularity.
\end{proof}

\providecommand{\bysame}{\leavevmode\hbox to3em{\hrulefill}\thinspace}

\end{document}